\documentstyle[12pt]{ article}
\begin{document}

\newtheorem{theorem}{Theorem}[section]
\newtheorem{prop}[theorem]{Theorem}
\newenvironment{mtheorem}{\begin{prop}\rm}{\end{prop}}
\newtheorem{coroll}[theorem]{Corollary}
\newenvironment{mcorollary}{\begin{coroll}\rm}{\end{coroll}}
\newtheorem{lemm}[theorem]{Lemma}
\newenvironment{mlemma}{\begin{lemm}\rm}{\end{lemm}}
\newtheorem{anoprop}[theorem]{Proposition}
\newenvironment{mproposition}{\begin{anoprop}\rm}{\end{anoprop}}
\newenvironment{mproof}{\begin{trivlist}\item[]{\em
Proof: }}{\hfill$\Box$\end{trivlist}}
\newenvironment{mdefinition}{\begin{trivlist}\item[]
{\em Definition.}}{\hfill$\Box$\end{trivlist}}
\newtheorem{eg}{\rm\sl \uppercase{Example}}[section]
\newenvironment{example}{\begin{eg}\rm}{\hfill$\Box$\end{eg}}

\newcommand{\ca}{{\cal A}}
\newcommand{\cb}{{\cal B}}
\newcommand{\cc}{{\cal C}}
\newcommand{\cd}{{\cal D}}
\newcommand{\cf}{{\cal F}}
\newcommand{\cg}{{\cal G}}
\newcommand{\ch}{{\cal H}}
\newcommand{\cl}{{\cal L}}
\newcommand{\cn}{{\cal N}}
\newcommand{\cm}{{\cal M}}
\newcommand{\cs}{{\cal S}}

\newcommand{\csa}{{$C^*$-algebra}}
\newcommand{\nsa}{{non-self-adjoint}}
\newcommand{\cssa}{{$C^*$-subalgebra}}

\def \IR{\hbox{{\rm I}\kern-.2em\hbox{{\rm R}}}}
\def \iR{\hbox{{\sevenrm I\kern-.2em\hbox{\sevenrm R}}}}
\def \IN{\hbox{{\rm I}\kern-.2em\hbox{\rm N}}}
\def \IC{\hbox{{\rm I}\kern-.6em\hbox{\bf C}}}
\def \IQ{\hbox{{\rm I}\kern-.6em\hbox{\bf Q}}}
\def \ZZ{\hbox{{\rm Z}\kern-.4em\hbox{\rm Z}}}

\newcommand{\ga}{{\alpha}}
\newcommand{\gb}{{\beta}}
\newcommand{\gd}{{\delta}}
\newcommand{\gA}{{\Alpha}}
\newcommand{\gB}{{\Beta}}
\newcommand{\gG}{{\Gamma}}
\newcommand{\gD}{{\Delta}}
\newcommand{\gs}{{\sigma}}
\newcommand{\gS}{{\Sigma}}

\newcommand{\ot}{\otimes}
\newcommand{\op}{\oplus}
\newcommand{\pr}{\prime}

\begin{center}
{\Large\bf
Homology for Operator Algebras II :
\vspace{.2in}

Stable Homology for Non-self-adjoint
\vspace{.15in}

Algebras}
\end{center}

\vspace{.3in}
\begin{center}
{\Large S. C. Power}

\vspace{.3in}

\it Department of Mathematics\\
Lancaster University\\
England LA1 4YF
\rm
\end{center}
\vspace{1cm}

\sf
\begin{abstract}
A new homology  is defined for a non-self-adjoint operator
algebra with a distinguished masa which is based upon
cycles and boundaries associated
with complexes of partial isometries in the stable algebra.
Under natural hypotheses the zeroth order group coincides
with the $\ K_0$\ group of the generated \csa . Several identifications
and applications are given, and in particular it is shown how
stable homology is significant for the classification
of regular subalgebras and regular limit algebras.
\end{abstract}
\rm

\newpage
Non-self-adjoint operator algebras are usually given in terms of a
construct from a more primitive category. Such categories include
partially ordered measure spaces (for nest algebras and commutative
subspace lattice algebras), semigroup actions (for semicrossed
products), ordered Bratteli diagrams (for subalgebras of AF
$\ C^*$-algebras), and binary relations and groupoids (for various
subalgebras of coordinatised von Neumann algebras and $\ C^*$-algebras).
Throughout the literature there has been a great emphasis placed on
relating operator algebras to the pertinent aspects of their genesis in
the simpler category. In the present paper we introduce various stable homology
groups $\ H_n ({\cal A}; \cc)$\  for operator algebras
$\ {\cal A}$\  with a prescribed self-adjoint subalgebra $\ \cc$.
The case of a digraph algebra provides the root context and here
stable homology is coincident with the integral
simplicial homology of the
simplicial complex of the underlying directed graph.
Although
intrinsically defined the stable homology groups, in contrast to
those of Hochschild cohomology, are often instantly
computable from the
underlying construct.  At the same time these groups are
related significantly to the algebraic structure.
On the other hand,
it is immediate from the definition that stable homology provides isomorphism
invariants for the most natural isomorphisms,
namely those  with
$\ C^*$-algebra extensions.

Although the new homology groups are of interest in their own right,
and in counterpoint with Hochschild cohomology,
they acquire added significance
with regard to the classification of
so-called {\em regular } subalgebras of non-self-adjoint operator algebras.
The results here are of interest even in the finite-dimensional case.
All self-adjoint subalgebras of finite-dimensional \csa s are
regular and they are well-understood
in terms of Bratteli diagrams.
The non-self-adjoint generalisation of this is to understand
regular subalgebras of digraph algebras, both in terms of
generalised Bratteli diagrams, and in terms of induced maps
on the  $\ K_0$ \  and homology groups,
together with other related invariants.
This is necessary to describe
not merely the nature of subalgebras,  but also their possible
positions, that is their classification up to inner
conjugacy. Even in the case of rather simple
digraphs, such as the cube $\ Cu$\ of section 3,
this raises some interesting combinatorial problems.

The structure and classification of regular subalgebras in finite dimensions
is also a necessary prelude to the classification of
limit algebras (even algebraic direct limits)
in the style of Elliott's classification
of AF \csa s in terms of the scaled $\ K_0$\  group.
Such ideas have already appeared in
\cite{scp-book} where it has been shown how
certain limit homology groups arise in the analysis of limits of
digraph algebras based on cycles.
The homology group formulations  below
give an alternative more generally applicable approach
to these limit groups.

The underlying idea for stable homology is simply the
following. Self-adjoint projections can play the role of
$0$-simplices, and partial isometries can play the role of
$1$-simplices. The formulation should provide a zeroth
order homology group that is coincident with the $\ K_0$\ group for the
generated \csa . (In the case of  a digraph algebra this is
a free abelian group whose rank is equal  to the number of components
of the digraph.) And the formulation should provide nonzero elements
in the first homology group if there are (appropriate) cycles
of partial isometries which are not expressible as boundaries in any
larger supercomplex. In this fashion we can obtain a homology
theory in which we can identify contributions from partial
isometry cycles that are linked to specific elements of (the positive
cone of) the $\ K_0$\ group of the generated \csa .
In brief, we define $\ H_n ({\cal A}; \cc)$\   in terms of the
simplicial homology of certain cycles of $\ \cc$-normalising
partial isometries in the stable algebra of
$\ {\cal A}$.
Although we concentrate on operator algebras in the
text this geometric form of homology is also
applicable to subspaces of \csa s which are
bimodules for a distinguished self-adjoint subalgebra.
In most of the examples we look at the
elements of $\ H_n ({\cal A}; \cc)$\  are already generated by
partial isometries of $\ {\cal A}$\ , rather than partial isometries of the
stable algebra. Nevertheless the stable algebra
formulation seems to be appropriate
for the purposes of classification
and of providing computable higher order
obstructions to the vanishing of Hochschild cohomology.
If $\ \ca$\  itself is self-adjoint then these groups vanish for
$\ n \ge 1$\ , {\em for entirely trivial reasons}. This already
suggests that these invariants are
particularly appropriate to general operator algebras and are intrinsically
more computable than Hochschild cohomology.

The paper is organised in the following way. In the first section we
define stable homology and remark on the
inadequacies
of some
variants of this definition. In section 2 we identify stable homology
in some fundamental settings
(i),  the tensor product of a digraph algebra
and a general \csa \ (Theorem 2.1 provides an elementary Kunneth formula), \
(ii) non-self-adjoint subalgebras
of factors determined by a finite lattice of commuting
projections,
and (iii)
regular limits of digraph algebras. In the latter case we recover
the homology limit groups
introduced in Davidson and Power \cite{dav-scp}.
We also mention a connection between the first stable homology group
and certain locally inner automorphisms.

In the remainder of the text we give three related applications.
Section 3 is concerned entirely with finite-dimensional
matters : regular subalgebras (and inclusions) of digraph algebras,
rigid embeddings, and the $\ K_0 \oplus H_*$\  uniqueness property,
particularly
in the context of
cycles, suspensions, discrete tori, and the cube algebra (this being a higher
dimensional variant of the 4-cycle algebra). In section 4
we indicate how such classifications
may be extended to similar settings in AF \csa s
by considering {\em scaled}  homology groups. In
the final section we
illustrate how homology can appear in
the classification
of regular limit algebras. It is clear that there are some very interesting
classification
problems in this area and we hope to develop these
ideas more fully elsewhere.
Note that the final two sections are independent of the first two in the sense
that one can consider the homology groups there as limit groups
(cf. Theorem 2.5).

Let us remark very briefly on the current literature concerning
the cohomology and homology of non-self-adjoint
operator algebras.

Automorphisms and derivations have formed a central
topic in operator algebra - one which is closely connected to
the more general considerations of Hochschild
cohomology. In the realm of reflexive algebras
the vanishing of Hochschild cohomolgy for nest algebras
has been demonstrated by Lance \cite{lan} and Christensen
\cite{chr}, whilst nonzero cohomology and non-inner derivations have been
identified and studied  by Gilfeather \cite{gil},
Gilfeather, Hopenwasser and Larson \cite{gil-hop-lar},
Gilfeather and Moore \cite{gil-moo}, and Power \cite{scp-spectral}.

Traditional studies of
Hochschild  cohomology for function algebras, as propounded by
Helemskii  \cite{hel}, Johnson \cite{joh} and Taylor \cite{tay},
for example, have direct bearing on operator algebras in the abelian
case. However a number of more recent studies have been pointed
specifically towards noncommutative algebras. In particular
Gilfeather and Smith \cite{gil-smi-1}, \cite{gil-smi-2} and
\cite{gil-smi-3} have examined Hochschild cohomology for constructions of
operator algebras analogous to the join, cone and
suspension constructions that are available in
simplicial homology. This work was inspired partly
by the cohomological identifications of
Gerstenhaber and Schack \cite{ger-sch} and Kraus and Schack \cite{kra-sch}
who promoted the fact that (for digraph algebras) Hochschild cohomology
is identifiable with  a simplicial cohomology.
The analysis of \cite{gil-smi-1}, \cite{gil-smi-2} and
\cite{gil-smi-3} also leans on basic techniques of
Johnson, Kadison and Ringrose (\cite{joh-kad-rin},
\cite{kad-rin}) in the Hochschild cohomology
theory for \csa s and von Neumann algebras.

In a different direction, but also motivated by digraph
algebras, Davidson and Power \cite{dav-scp}
considered direct limit homology groups
for triangular limit algebras, and in \cite{scp-book} it was shown that
these could be used as classifying invariants in certain contexts of
4-cycle limit algebras. This homology theory, like that which we have
given for reflexive algebras in \cite{scp-spectral},
is closely tied to the underlying
coordinatisation of the algebra, and is
possibly more appropriate and computable than Banach
algebra cohomology.
The present paper gives a general intrinsic
formulation for these groups
which is quite widely applicable.
We envisage that these invariants will be accessible and significant
in the area of
subalgebras of groupoid \csa s, as developed by Muhly
and Solel \cite{muh-sol} in the triangular case, and in
the (completely undeveloped)
area of direct limits of non-self-adjoint
subhomogeneous algebras.

The following terminology is adopted.
A digraph algebra $\ \ca$\  is a subalgebra of
a complex matrix algebra $\ M_n$\
which contains a maximal abelian self-adjoint subalgebra (a masa).
These are also known as finite-dimensional CSL algebras
or finite-dimensional incidence algebras.
If $\ \{e_{i,j}\}$\  is a standard matrix unit  system for $\ M_n$\
such that the masa in question is spanned by the matrix units $\ \{e_{i,i}\}$\
then the digraph for $\ \ca$\   has $\ n$\  vertices and directed
edges $\ (i,j)$\  for each $\ e_{i,j}$\  in $\ \ca$. \  This
digraph (or binary relation) is transitive and reflexive,
with no multiple directed
edges. From the point of view of cohomology and homology the digraph
algebras $\ A(D_{2n})$\  for the $\ 2n$-cycle digraphs $\ D_{2n}$\  are the
first significant examples.
These algebra are also denoted $\ \ca_{2n}$\  and are occasionally refered to
as the (finite-dimensional)
tridiagonal matrix algebras. All the algebras that we
consider are viewed as subalgebras
of \csa s, and by a star-extendible homomorphism we mean one which
is a restriction of a
\csa \ homomorphism between the generated \csa s.

\newpage

\section{Formulations of Stable Homology}

Let $\ \ca$\   be an operator algebra with self-adjoint subalgebra
$\ \cc$.\  Our main interest is when $\ \cc$\  is maximal abelian.
In the following discussion it should be held in mind that we seek to
formulate a stable homology theory, in the sense that
$\ H_n(\ca \otimes M_n; \cc \otimes \ \IC^n) = H_n(\ca ; \cc)$,\
we wish to have $\ H_0(\ca ; \cc) = K_0(C^*(A)$\ ) in appropriate contexts,
and we require that $\ H_*(\ca;\cc)$\  specialises to integral simplicial
homology in the case of digraph algebras.
Furthermore, we wish to have the elementary Kunneth formula of
Theorem 2.1 which links $\ K_0$\  and $\ H_*$.\

An alternative elementary formulation of $\ H_1(\ca ; \cc)$,\  which is
independent of simplicial homology, is given in Remark 1.3.

The stable algebra of an operator algebra $\ \ca$\ is taken to be the
algebra of  finitely nonzero  infinite matrices over $\ \ca$.\
Let $\ B$\  be a
finite-dimensional $\ C^*$\ -algebra contained in the stable algebra
$\ M_\infty (C^*(\ca))$\  with a full matrix unit system $\ \{f_{ij}\}$\
consisting of
partial isometries which normalise the
subalgebra  $\ M_\infty
(\cc)$.\
This means that if $\ f \in \{f_{i,j}\} $\ and $ c \in \ M_\infty
(\cc)$\  then $\ fcf^*$\ and
$\ f^*cf$\ belong to $\ M_\infty
(\cc)$.\
Then the subalgebra

\[
A = B \cap M_\infty ({\ca})
\]

\noindent contains the diagonal matrix units and so $\ A$\
is a subalgebra of $\ B$\
associated with the binary relation $\ R(A) = \{(i, j): f_{i,j} \in
A\}$.\
In particular $\ A$\  is completely isometrically isomorphic to the digraph
algebra
associated with $\ R(A)$.\
If $\ A$\  and $\ A^\prime$\  are two such subalgebras of $\ M_\infty (\ca)$\
then we declare them to be equivalent if they are conjugate
by means of a unitary operator $\ u$\  in (the unitisation of)
$\ M_\infty ({\ca} \cap {\ca}^*)$.\

Let
[A] denote the equivalence class of such digraph
subalgebras, and let $\ H_n ([A])$\
denote the $\ n^{th}$\  integral simplical homology group of the simplicial
complex $\ \Delta ([A])$\  associated with $\ R(A)$.\
The complex $\ \Delta ([A])$\  is perhaps most easily specified by
viewing $\ R(A)$\  as the edges of a directed graph $\ G$\  with
vertices $\ v_1,...,v_n$\  :
Let $\ \overline{G}$\  be the undirected graph of $\ G$.\  Then the
0-simplices of $\ \Delta(G)$,\  denoted
$\ \sigma_i =\  <v_i>$\  , $\ 1 \leq i \leq n$,\  are associated
with the vertices $\ v_i$\  of $\ G$,\  and the t-simplices of $\ \Delta(G)$\
correspond to
the complete subgraphs of $\ \overline{G}$\  with $\ t + 1$\  vertices.
Thus if $\ v_i, v_j, v_k$\  determine a
complete subgraph of $\ \overline{G}$\  then the 2-simplex $\ \sigma_{ijk} =
<v_i,
v_j, v_k>$\  is included in $\ \Delta(G)$.\

The group  $\ H_n
({\ca} ; \cc)$\  is defined to be the quotient

\[
H_n ({A} ; \cc) = (\sum_{[A]} \oplus H_n ([A]))/Q_n
\]

\noindent where the direct sum indicates the
restricted direct sum, and where $\ Q_n$\
is a natural subgroup corresponding to inclusion identifications and to
identifications arising from certain orthogonal direct sums
(induced decompositions) as described below.
Roughly speaking, it follows that $\ H_1 ({\ca}; \cc)$\  is nonzero
if there exists a sequence of normalising
partial isometries in
$\ M_\infty ({\ca})$\  which form a 1-cycle in a finite-dimensional algebra
$\ A$\  but which do not give a 1-boundary in any affiliated algebra
$\ A^\prime$\  containing $\ A$.\

We now define $\ Q_n$.\  Refer to the algebras $\ A, A^\prime$,\  as above, as
$M_\infty ({\cc})$-{\it normalising} (or \ $\cc$-{\it normalising})
{\it digraph} {\it algebras}
{\it for} $\ {\ca}$,\
and refer to the matrix
unit system $\ \{f_{i,j}: f_{i,j} \in A\}$\  as a {\it  partial}
{\it matrix} {\it unit}
{\it system for} $ \ A$.\  Note that such a  system has
the special property that the
generated star semigroup is a full matrix unit system in the usual
sense. Let $\ A \subseteq  A^\prime$\  \ be  $\ \cc$\ -normalising digraph
algebras
such that the partial matrix unit system of $\ A$\  is a subset of the
partial matrix unit system for $\ A^\prime$.\  Then there is a natural
well-defined group
homomorphism $\ \theta: H_n ([A]) \rightarrow H_n ([A^\prime])$\  which is
induced by the resulting digraph {\bf inclusion} $\ R(A) \rightarrow
R(A^\prime)$.\
Identify each group $\ H_n ([A])$\  with its summand in $\ \sum_{[A]} \oplus
H_n
([A])$\  and let $\ Q^a_n$\  be the set of elements of the form \ \ $\ g -
\theta
(g)$\  \ \ associated with all such group homomorphisms $\ \theta : H_n ([A])
\rightarrow H_n ([A^\prime])$,\  and elements $\ g$\  in $\ H_n ([A])$.\
Of course there may be a finite number of such group homomorphisms for each
pair
$\ [A], [A^\prime]$.\
Note, in particular, that we only consider rather special inclusions
which, in the terminology of section 3, are multiplicity one regular
inclusions.

The
subgroup $\ Q_n$\  is defined to be the subgroup
generated by $\ Q^a_n$\  and $\ Q^b_n$,\  where $\ Q_n^b$\
corresponds to certain orthogonal direct sum identifications, as we now
indicate.

Let $\ A$\  be a $\ \cc$-normalising digraph algebra for  $\ \ca$\  with
partial matrix unit
system $\ \{f_{i,j}: (i,j) \in R(A)\}$.\  Without loss of generality assume
that $\ C^* (A) = M_n$.\  Let $\ f_{11} = f^\prime_{11} \oplus f^\prime_{22}$,\
with $\ f^\prime_{11}, f^\prime_{22}$\  nonzero projections in
$\ M_\infty (\cc)$.\  Then,
since the $\ f_{i,j}$\ are  $\ \cc$-normalising, it follows that
there is an {\bf induced decomposition}
$\ f_{ij} = f^\prime_{ij} + f^{\prime \prime}_{ij}$,\  for $\ (i,j)$\  in
$\ R(A)$,\  such that

\[
\ \{f^\prime_{ij} : (i,j) \in R(A)\}\  \ \ \ \mbox{and} \ \ \  \ \{f^{\prime
\prime}_{ij} : (i,j) \in R(A)\}
\]

\noindent are partial matrix unit systems for
$\ \cc$-normalising digraph algebras $\ A^\prime, A^{\prime \prime}$\
respectively.
In fact $\ f^\prime_{ij} = f_{i,1}f^\prime_{1,1}f_{1,j}.$\
Let
$\ \theta^\prime: A \rightarrow A^\prime, \ \ \theta^{\prime \prime} :
A \rightarrow
A^{\prime \prime}$\  be the associated algebra isomorphisms, with induced
(well-defined)
isomomorphisms

\[
\theta^\prime_n : H_n ([A]) \rightarrow  H_n ([A^\prime]),
\ \   \theta_n^{\prime \prime}: H_n ([A]) \rightarrow H_n
([A^{\prime \prime}]).
\]

\noindent Define $\ Q^b_n$\  to be the set of elements of
$\ \sum_{[A]} \oplus H_n ([A])$\  of the form

\[
g - \theta^\prime_n (g) - \theta^{\prime \prime}_n (g) \ \ \ \ , g \in
H_n([A]).
\]

The definition of the stable homology groups is now complete.

\vspace{.3in}
\noindent {\bf Definition 1.1} \ The $\ \cc$\ -{\it normal stable homology} of
the
operator algebra $\ {\cal A}$\  with distinguished self-adjoint
subalgebra $\ \cc$\  consists of the groups $\ H_n ({\cal A} ; \cc)$,\  $\ n =
0, 1, 2, \dots$.\
\vspace{.3in}

The discussion above gives a fairly intuitive
construction  and we shall see in the next section that it is quite suited to
calculations in specific contexts.

Let us point out how the Grothendieck group $\ G(S)$\  \ of an abelian
unital semigroup $\ S$\  can be viewed in the above formalism.
Let
\[
{\cal G} = (\sum_{s\in S} \op \ZZ)\slash {\cal R},
\]
where ${\cal R}$\
is the subgroup generated by the elements associated with the relations for the
semigroup $\ S.$\ (A typical such element has the form
$\  n_{s+t} - (n_s \op n_t)$\  with $\ s, t $ in $\ S$.)\
Then
$\ \cg$\  is naturally
isomorphic to the usual Grothendieck group of $\ S$.\
{}From this and our definitions above
it follows that if $\ \cb$\  is a unital \csa \ then $\ H_0(\cb;\cb) =
K_0(\cb).$\

Similarly,
let $\ \cc$\  be a unital \ \cssa \ of $\ \cb$\  with the following properties
:
(i) for each projection class
$\ [p]$\ in $\ M_\infty(\cb)$\ there is a projection $\ q$\ in
$M_\infty(\cc)$\ with $\ [p] = [q]$,
and (ii) if $\ q_1$\  and $\ q_2$\ are projections in \ $M_\infty(\cc)$\
which are equivalent in \ $M_\infty(B)$\ then they are equivalent
by an $M_\infty(\cc)$-normalising partial isometry. The first property implies
that
the natural map $\ K_0(\cc) \to K_0(\cb)$\  is a surjection.
If (i) and (ii) both hold we shall say that the map
$\ K_0(\cc) \to K_0(\cb)$\  is a {\em regular surjection}.
Under these circumstances it follows that
$\ H_0(\cb;\cc) = K_0(\cb).$
\vspace{.3in}

\noindent {\bf Remark 1.2}\ \
One can also
present the homology groups $\ H_n (\ca ; \cc)$\  in a more orthodox fashion
as the homology groups of a chain complex $\ (S_n ({\cal A}), d_n)$.\
To do this define $\ S_n ({\cal A})$\  to be the quotient

\[
(\sum_{[A]} \oplus S_n ([A]))/QS_n
\]

\noindent where $\ S_n ([A])$\  is the $\ n$\ -chain group of the complex for
$\ R(A)$,\  with
integral coefficients, and where $\ QS_n$\  is the subgroup determined by
the relations of inclusion of matrix unit systems and of orthogonal
direct sum. The boundary operators $\ d_n$\
respect the subgroups $\ QS_n$\ and so we may define
$\ Z_n ({\cal A})$,\  the
$\ n$\ -cycle group, and $\ B_n ({\cal A})$,\  the $\ n$\ -boundary group. Then
the quotient groups $\ Z_n
({\cal A})/B_n ({\cal A})$\  are the homology groups of the associated
chain complex $\ (S_n ({\cal A}), d_n)$\  and they are identifiable
with the groups $\ H_n ({\cal A}, \cc)$.\
\vspace{.3in}

\noindent {\bf Remark 1.3}
\ \ An alternative direct formulation of $\ H_1(\ca ; \cc)$\
can be made in the following fashion.

A {\it basic $\ M_{\infty}(\cc)$-normalising 1-cycle} of $\ \ca$\  is a triple
$\ \sigma = (u_1,u_2,u_3)$,\  or a $\ 2n$\ -tuple $\ \sigma =
(u_1,\dots,u_{2n})$,\
consisting of partial isometries in $\ M_{\infty}(\ca)$\
which normalise $\ M_{\infty}(\cc)$\  and satisfy the relations
suggested by the following diagrams.
\vspace{.3in}

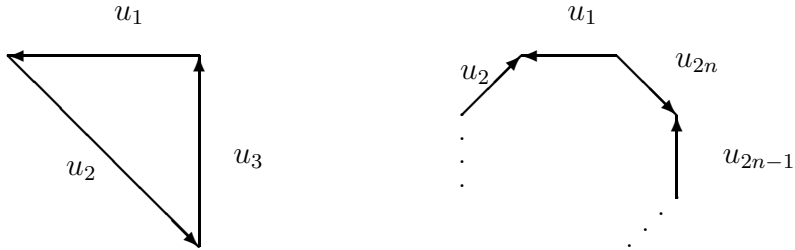
\setlength{\unitlength}{0.0125in}%
\begin{picture}(90,104)(225,680)
\thicklines
\put(305,680){\vector( 0, 1){ 80}}
\put(225,760){\vector( 1,-1){ 80}}
\put(265,775){\makebox(0,0)[lb]{\raisebox{0pt}[0pt][0pt]{\twlrm $\ u_1$\ }}}
\put(245,710){\makebox(0,0)[lb]{\raisebox{0pt}[0pt][0pt]{\twlrm $\ u_2$\ }}}
\put(315,715){\makebox(0,0)[lb]{\raisebox{0pt}[0pt][0pt]{\twlrm $\ u_3$\ }}}
\put(480,760){\vector(-1, 0){ 40}}
\put(480,760){\vector( 1,-1){ 25}}
\put(305,760){\vector(-1, 0){ 80}}
\put(505,700){\vector( 0, 1){ 35}}
\put(500,755){\makebox(0,0)[lb]{\raisebox{0pt}[0pt][0pt]{\twlrm $\ u_{2n}$\ }}}
\put(415,735){\vector( 1, 1){ 25}}
\multiput(415,735)(0.00000,-10.00000){4}{\makebox(0.4444,0.6667){\tenrm .}}
\multiput(505,700)(-6.66667,-6.66667){4}{\makebox(0.4444,0.6667){\tenrm .}}
\put(455,775){\makebox(0,0)[lb]{\raisebox{0pt}[0pt][0pt]{\twlrm $\ u_1$\ }}}
\put(410,750){\makebox(0,0)[lb]{\raisebox{0pt}[0pt][0pt]{\twlrm $\ u_2$\ }}}
\put(520,715){\makebox(0,0)[lb]{\raisebox{0pt}[0pt][0pt]{\twlrm $\ u_{2n-1}$\
}}}
\end{picture}

\vspace{.3in}

Thus, for the $\ 2n$\ -cycle, $\ d(u_{2k}) = d(u_{2k+1})$\
and $\ r(u_{2k+1}) = r(u_{2k+2})$\  for all appropriate $\ k$,\  and all these
domain and range projections are orthogonal. Furthermore,
\[
u_{2n} = u_{2n-1}u_{2n-2}^* \dots u_{2}^*u_{1},
\]
and so $\ C^*(\{u_1,\dots,u_{2n}\})$\  is isomorphic to $\ M_{2n}$,\
and the nonzero words in the elements
$\ u_1,\dots,u_{2n}$\  and their adjoints provide a complete matrix unit
system for the algebra.

Let $\ [\sigma]$\  denote the set of basic $\ M_{\infty}(\cc)$\ -normalising
$\ 1-$\ cycles \ \ $\ \omega$\  \ \ that are
unitarily equivalent to $\ \sigma$\  in the sense
that $\ \omega = z\sigma z^*$\
for some unitary $\ z$\  in (the unitisation of)
$\ M_{\infty}(\ca) \cap M_{\infty}(\ca)^*$.\
Define $\ Z_1(\ca ; \cc)$\  to be the free abelian
group generated by such classes, modulo the following three relations:

(i) (orthogonal sum)
\[
[(u_1,\dots,u_{2n})] + [(v_1,\dots,v_{2n})] =
[(u_1+v_1,\dots,u_{2n}+v_{2n})],
\]
where representatives are chosen so that $\ u_i + v_i$\  is a
partial isometry for all $\ i$.\

(ii) (cancellation)
\[
[(u_1,\dots,u_{2n})] + [(u_{2n},\dots,u_{1})] = 0.
\]

(iii) (addition)
\ \ If
$\ \sigma_1 =
(u_1,\dots,u_{2n}),\ \  \sigma_2 = (v_1,\dots,v_{2m}), \ \ u_{2n} = v_1,
$\  and$\  \ \ \sigma = (u_1,\dots,u_{2n-1},v_2,\dots,v_{2m}),$\
then
\[
[\sigma] = [\sigma_1] + [\sigma_2].
\]

Define $\ B_1(\ca ; \cc)$\  to be the subgroup generated by the classes of
the $\ 1$-cycles $\ \sigma$\  coming from triples. Then
$\ H_1(\ca ; \cc) = Z_1(\ca ; \cc) \slash B_1(\ca ; \cc)$.\

\vspace{.3in}

\noindent {\bf Remark 1.4} \ It is tempting to drop the normalising
condition in the above formulations
and define a stable homology in terms of
all digraph subalgebras of $\ M_{\infty}(\ca \bigcap \ca^*)$\  with respect
to unitary equivalence from $\ M_{\infty}(\ca \bigcap \ca^*)$.\
But this move leads to unwanted complications
in view of the proliferation of unitary equivalence classes of
partial isometry cycles, even
when $\ \ca$\ is a digraph algebra.
In fact one does not obtain a homology theory which generalises simplicial
homology
in this case. For example, in the case of the basic digraph algebra
$\ \ca = A(D_4)$\ \ \  the resulting $\ H_1$\  group  is
the restricted direct product of uncountably many copies of $\ \ZZ$.\
This is related to the fact that there are uncountably many
inner equivalence classes of partial isometries in this algebra.
\vspace{.3in}

\noindent {\bf Remark 1.5}\  \ Here are three variations of stable homology:

(i) One could be more restrictive in the choice of partial matrix units
by demanding that they normalise the diagonal algebra
$\ D_\infty(C)$\  rather than
$\ M_\infty(C)$.\
This homology is somewhat more computable
and it is adequate
for the approximately finite settings considered
in section 3, 4 and 5.
However, there is the big disadvantage that one
does not obtain a variant of Theorem 2.1
below.

(ii) One could drop the dependence on $\ \cc$\  altogether and define the
homology groups $\ H_n(\ca;\ca \cap \ca^*)$.\
Here one requires normalisation of $\ M_{\infty}(\ca \cap \ca^*).$\
This is a very attractive move, superficially,
since the resulting groups are invariants for
star-extendible isomorphism.
Furthermore this homology does coincide with the
simplicial homology of the digraphs of the digraph algebras.
(See Theorem 2.2 (i).) However, the functoriality properties are seriously
inadequate in the sense that
regular morphisms between digraph algebras (such as the
rigid embeddings in section 3)
do not induce homology group homomorphisms.
Furthermore in many basic contexts of interest these homology groups
are clearly inappropriate.
To see this consider the following example.

Ler $\ \cf $\ be the direct limit algebra
$\displaystyle{\lim_\to(M_{2^k},\phi_k)}$\  (not necessarily closed)
where $\ \phi_k(a) = a \op a$\  for all $\ a$.\
Let $\ \ca $\  be the subalgebra of $\ A(D_4) \ot \cf $\  consisting of the
operators
$\ a$\  for which $\ (e_{1,1} \ot 1)a(e_{1,1} \ot 1$) \ belongs to
$\ e_{1,1} \ot \cd $\ where
$\  \cd =
\displaystyle{\lim_\to(D_{2^k},\phi_k)},$\  the standard diagonal subalgebra.
For the natural masa $\ \cc = \IC^4 \ot \cd$\  the normal stable homology
$\ H_1(\ca;\cc)$\ is nontrivial and can be readily identified
using Theorem 2.5. On the other hand $\ H_1(\ca;\ca \cap \ca^*)$\ is trivial.
This is essentially because the normalising demand is too great; if
$\ v$\  is a partial isometry which normalises $\ \ca \cap \ca^*$\
then
$\ (e_{1,1} \ot \ 1)v(e_{3,3} \ot \ 1) =
\ (e_{1,1} \ot \ 1)v(e_{4,4} \ot \ 1)  = 0.$

(iii) One could {\em restrict} the class of partial
isometries that are admisssible in the partial matrix
unit systems of the digraph subalgebras. For example, in the
operator algebra
of Example 2.3 restriction to finite rank matrix units leads
to a trivial first restricted stable homology
group, and this reflects the triviality of the first
simplicial homology group
of the associated digraph of that example. This type
of restriction
seems  appropriate for an analysis of the homology
affiliated to elements
of $\ K_0(C^*(\ca))$.\
\vspace{.3in}

\noindent {\bf Remark 1.6} \ \ The stable homology that we have given is
defined in terms of finite-dimensional
$\ C^*$-subalgebras. Even in the case
triangular limit algebras
with an "approximately finite-dimensional character"
such an "AF
homology" may
be inappropriate. We have in mind here
the limits of cycle algebras
under {\em non-star-extendible} embeddings,
given in \cite{scp-simplicial} and \cite{scp-book}.
It can be shown that these have trivial first stable homology
(with resect to the unique masa). On the other
hand they do posses natural nonzero
limit homology groups (see \cite{scp-simplicial}).
\vspace{.3in}

\noindent {\bf Remark 1.7}
Minor modifications of the definitions above lead to the formulation
of the relative homology groups :

Let $\ \cc, \cal A$\  be before, and let $\ \cal A^\prime$\  be an
intermediate operator algebra with $\ \cc \subseteq {\cal A^\prime} \subseteq
{\cal A}$.\  Let $\ B, A$\  be as before, and let $\ A^\prime = B \cap M_\infty
({\cal A}^\prime)$,\  so that $\ A^\prime$\  is a $\ \cc$-normalising digraph
subalgebra of $\ {\cal A}^\prime$\  which is spanned by some of the matrix
units of $\ A$.\  To the unitary equivalence class $\ [A, A^\prime]$\  of such
pairs associate the relative integral simplicial homology group
$\ H_n ([A, A^\prime])$,\  which is defined to be the relative homology
group $\ H_n (\Delta (A), \Delta (A^\prime))$,\  where $\ \Delta (A^\prime)$\
is the
subcomplex determined by $\ R(A^\prime)$.\  Define the relative
$\ \cc$-normalising homology to be the quotient

\[
H_n ({\cal A}, {\cal A}^\prime; \cc) = (\sum_{[A, A^\prime]} \oplus H_n
([A, A^\prime]))/Q_n ({\cal A}, {\cal A}^\prime)
\]

\noindent where $\ Q_n({\cal A}, {\cal A}^\prime)$\  is the subgroup of the
restricted
direct sum determined by orthogonal direct sum identifications, and by
subcomplex identifications.

Alternatively, we can view the chain complex $\ (S_n ({\cal A}^\prime),
d_n)$\  as a subcomplex of the chain complex $\ (S_n ({\cal A}), d_n)$,\  in
which case $\ H_n ({\cal A}, {\cal A}^\prime; \cc)$\  is the homology
of the quotient chain complex $\ (S_n ({\cal A})/S_n ({\cal A}^\prime),
d_n)$.\
\vspace{.3in}

\noindent {\bf Remark 1.8}\ \ Stable homology is, prima facie, an invariant
for pairs $\ (\ca,\cc)$.\  However, in the presence of uniqueness
theorems (up to automorphisms of $\ \ca$)\ for regular masas $\ \cc$,\ one can
simply define $\ H_{*}(\ca) =
\ H_{*}(\ca;\cc)$\  and obtain a well-defined
homology theory for $\ \ca$\  itself.
Examples of this appear in sections 3 and 4 and we expect similar definitions
of  $\ H_*(\ca)$\ in
much more general circumstances.
Of course, in the extreme case of triangular algebras, such
as the lexicographic products in
Example 2.6, the masa $\ \cc =  \ca \cap \ca^*$\ is intrinsic to the algebra
and we may define $\ H_*(\ca) = \ H_*(\ca;\cc)$.

\newpage

\section {\bf Identifications of Stable Homology}

We have remarked in the introduction that the stable homology
of a digraph algebra coincides with the simplicial
homology of the complex for the digraph of the algebra.
The next two theorems establish this and
give different more general versions of
this correspondence.
The proofs are essentially elementary and depend
on the decomposition of an arbitrary $\ M_\infty(\cc)$-normalising
digraph algebra in the stable algebra into an "parallel sum" of
ones that are unitarily equivalent to certain  easily visible
elementary digraph subalgebras.
\vspace{.3in}

\noindent {\bf Theorem 2.1} \ \ \it  Let $\ A(G)$\  be a digraph algebra
and let $\ \cb$\  be a unital
\csa \ with abelian unital self-adjoint subalgebra $\ C$\  such that
the inclusion $\  C \to \cb$\  induces a regular
surjection $\ K_0C \to K_0\cb$.\  Then,
for each $\ n \ge 0$,\
\[
H_n ((A(G) \otimes \cb;\  \IC^{|G|} \otimes C)) \ =
\ H_n (\Delta (G)) \otimes_{\ZZ} K_0(\cb).
\]
\rm

\begin{mproof} Let $\  \ca =  A(G) \otimes \cb,\   \cc = \IC^{|G|} \otimes C$.\
Here $\ \IC^{|G|}$\  is the diagonal subalgebra of
$\ A(G)$\ with respect to a fixed matrix unit
system $\{e_{i,j}: (i,j) \in E(G)\}$.\
We may assume that $\ G$\  is connected.
The main step is to reduce the quotient
expression for $\ H_n(\ca;\cc)$\  to one
involving a direct sum over standard type digraph subalgebras of the
form $\ A(G) \otimes q$\  where $\ q$\  is a projection in $\ M_N (C).$\

Let $\ A \subseteq M_N(A(G) \otimes \cb) = A(G) \otimes M_N (\cb)$\
be a digraph subalgebra with partial matrix
unit system $\ \{f_{k,l}\}$\  each element of
which normalises \ $M_N (\cc) = \IC^{|G|} \otimes M_N (C)$.\
Without loss of generality
assume that the digraph for $\ A$\  is connected
and that the full system of $\ \{f_{k,l}\}$\  is
$\{f_{k,l}: 1 \le k, l \le K\}.$

Note the following principle : if a $\ 2 \times 2 $\ operator matrix $\ v$\ is
a
partial isometry, say

\[
v = \left[\begin{array}{cc}
a&b\\c&d\end{array}\right],
\]

\noindent and if $\ vxv^*$\  is block diagonal when $\ x$\ is

\[
\left[\begin{array}{cc}
I_1&0\\0&0\end{array}\right]
\ \ \mbox{ and } \ \
\left[\begin{array}{cc}
0&0\\0&I_2\end{array}\right],
\]

\noindent then $\ a, b, c, d $\  are partial isometries with orthogonal
domains and ranges. Using this
principle repeatedly obtain an induced decomposition
$\ f_{k,l} = f_{k,l}^{(1)} + \dots + f_{k,l}^{(t)},$\
in the sense given in section 1,
such that each  $\ f_{k,l}^{(r)}$\  has the normalising property
and  belongs to one
of the spaces $\ e_{i,j} \otimes M_N (\cb)$.\
More explicitly, consider the projections $\ e_i = e_{i,i} \ot 1, \ 1 \le i \le
|G|$,\  in $\ A(G) \ot M_N(C)$.\  Then there is an induced decomposition
$\ f_{k,l} = f_{k,l}^\pr + f_{k,l}^{\pr\pr},$
where, for each pair $\ k, l,$\ \
$\ f_{k,l}^\pr = f_{k,1}(f_{1,1}e_1)f_{1,l}$.\
If this is a nontrivial decomposition, that is, if $\ f_{1,1}e_1 \ne 0,$\
then $\ f_{1,1}^\pr e_1 = f_{1,1}^\pr.$\
Furthermore, the systems $\{f_{k,l}^\pr\}$\ and $\{f_{k,l}^{\pr\pr} \}$\
still have the normalising property.
Repeating such decompositions leads to the desired reduction.

For  fixed $\ r$\  consider the associated {\em full}
matrix unit system  $\ \{f_{k,l}^{(r)}\}$.\
Then for each $\  i$\ the intersection
\[
\ \{f_{k,l}^{(r)}\} \cap ( e_{i,i} \otimes M_N (B))
\]
is a complete  system of matrix units and so has the form
\[
e_{i,i} \otimes g^i_{s,t}\ \ \ , \mbox{for}\ \ \  1 \le s, t \le n_i,
\]
where $\ \{g^i_{s,t}\}$\ is a full matrix unit system in $\ M_N (\cb)$\
which normalises $\ M_N(C)$.\  Let $\ m$\ be the maximum of
the  numbers $\ n_i$, say $\ m = n_p.$\
Since $\ f_{k,l}^{(r)}$\ is a full matrix unit system
the matrix unit $\ g^i_{s,s}$\  is equivalent to $\ g^p_{s,s}$\
for each $\ s$\ with $ \ 1 \le s \le n_i$,\  by a matrix unit of the form
\[
f_{k,l}^{(r)} = e_{i,p} \ot v
\]
where $\ v$\  normalises $\ M_N(\cc)$. It follows that by conjugating
we may assume that $\ g^i_{s,s} = \ g^p_{s,s}$\ for $\  1 \le s \le n_i.$\
We now see that the full matrix unit system $\ \{f_{k,l}^{(r)}\}$\
is conjugate,
by a normalising unitary in $\ M_N(\ca \cap \ca^*),$\
to a subsystem of a system of the form
\[
\{e_{i,j} \ot g_{s,t}\}
\]
where
$\ \{g_{s,t}\}$\ is a complete matrix unit system for $\ M_m$\
as a (not necessarily unital) subalgebra of $\ M_N (\cb)$,\
with the normalising property.

To recap, it has been shown that $\ A$\  is inner equivalent,
by a unitary in $\ M_N(\ca \cap \ca^*)$,\  to a digraph algebra
with partial matrix unit system
$\{f_{k,l}\}$\   admitting an induced decomposition
$\ f_{k,l} = f_{k,l}^{(1)} + \dots + f_{k,l}^{(t)},$\
where each partial system  $\{\  f_{k,l}^{(r)}\} $\  is a subsystem of the
standard system
for an elementary digraph algebra
$\ A(G) \ot M_m \ot q$,\  where $\ q$\  is a projection in $\ M_N(\cc)$,\
and $\ m, N $ and $ q$\  depend on $\ r$.\  In brief, each digraph
algebra class $\ [A]$\ for $\ \ca$\  has a representative digraph algebra
which is constructed in a natural way from elementary ones.

Let

\[
{\cal G}\ = \sum_{[A]} \oplus H_n([A]),\ \ \ {\cal G}_0 =
\sum_{[q],m} \oplus H_n([A(G) \otimes M_m \otimes q]).
\]

\noindent where $\ {\cal G}_0$\ is the subgroup of $\ \cg $\
associated with the elementary digraph
subalgebras indexed by the $\ K_0(\cb)$\ classes $\ [q]$,\  with $\ q$\  in
$\ M_\infty(C)$,\ and  positive integers $\ m$.\
Thus $\ H_n(\ca;\cc) = \cg/ Q_n$,\ and, by the reductions above,
$\  \cg/ Q_n =  \cg_0/ Q_n.$\
Furthermore,
$\  \cg_0/ Q_n =  \cg_0/ Q_{n,0} $\
where $\ Q_{n,0}$\  \ is the subgroup generated by the set
of relations $\  Q_{n,0}^a,\   Q_{n,0}^b$\   corresponding to
inclusions and induced decompositions for elementary digraph algebras.
This is purely algebraic fact which follows from the simple
principle that for abelian groups
$G, H$\ the quotient group
$\ (G \op G \op H)/\{g \op -g \op 0\}$\
is isomorphic to $\ 0 \op G \op H$.

The inclusion $\ A(G) \ot e_{1,1} \to A(G) \ot M_n$\  induces an
isomorphism of simplicial homology leading to the further reduction

\[
H_n(\ca;\cc) = (\sum_{[q]} \oplus H_n([A(G) \otimes q]))/Q_{n,0}^b
\]

\noindent where the direct sum extends over classes of projections
$\ q$\  in $\ M_\infty(\cc)$.
(There are no remaining inclusion relations.)
Thus, making the natural identifications
$\ H_n([A(G) \otimes q]) = H_n(\gD(G)),$\   we see that

\[
H_n(\ca;\cc) = (\sum_{[q]} \oplus H_n(\gD(G)))/S
\]

\noindent where $\ S$\  is the subgroup corresponding
to the semigroup relations for the
classes $\ [q]$.\
Hence

\[
(\sum_{[q]} \oplus H_n(\gD(G)))/S = H_n(\gD(G)) \ot_{\ZZ} ((\sum_{[q]} \op
\ZZ)/S).
\]

\noindent Since the map $\ K_0C \to K_0\cb$\  is a
regular inclusion it follows that

\[
K_0\cb = ( \sum_{[q]} \op \ZZ)/S
\]

\noindent and the proof is complete.
\end{mproof}
\vspace{.3in}

The next identifications are similar to the last
but are somewhat more elementary.

Let $\ \cm$\  be a factor and let $\ \cal L$\  be a finite lattice of commuting
projections in $\ \cm$\  with associated subalgebra $\ \cal A$\  consisting of
the
operators $\ a$\  in $\ \cm$\  for which $\ (1 - p) a p = 0$\  for all $\ p$\
in $\ \cal
L$.\  The minimal nonzero interval projections $\ { f - e}$,\  with $\  f
> e$\  projections of $\ \cal L$,\  form a finite set, $\ Q = \{q_1, \dots,
q_n\}$\  say. $\ Q$\  carries the transitive partial order $\ \ll$\  where

\[
q \ll q^\prime \ \ \ \Leftrightarrow \ \ \ q {\cal A} q^\prime = q\cm q^\prime.
\]

\noindent Write $\ H_n (\Delta({\cal L}))$\  for the integral simplicial
homology
of the complex $\ \Delta({\cal L})$\  for the partial order $\ \ll$,\  viewed
as a digraph.
\vspace{.3in}

\noindent {\bf Theorem 2.2} \it Let $\ \cm$\  be II$\ _1$\  factor
and let $\ {\cal L} \subseteq
{\cal M}$\  be a finite lattice of commuting projections with
associated reflexive algebra
$\ \cal{A} \subseteq
{\cal M}$.\  Then

\noindent (i)
\[
H_n ({\cal A; A \cap A^*}) =
H_n (\Delta({\cal L}))
\otimes_{\ZZ} \IR.
\]

\noindent (ii) If $\ {\cal C} \subseteq \cm$\  is a regular masa of $\ \cm$\
then
\[
H_n ({\cal A}; \cc) = H_n (\Delta({\cal L}))
\otimes_{\ZZ} \IR.
\]
\rm

\begin{mproof}
(i)
Let $\ A$\ be a digraph algebra for $\ \ca$\ which is contained in
$\ M_N(\ca) = \ca \ot M_N$\
and has a partial matrix unit system
$\ \{f_{k,l}\}$\  which is elementary in the sense that for each
pair $\ k, l$\  the operator $\ (q_i \ot I_N)f_{k,l}(q_j \ot I_N)$\
is nonzero for at most one pair $\ i, j.$\  The conjugacy class
of each such subalgebra
is determined by a subdigraph $\ H$\ of $\ G$\ and a projection
$\ q$\ in $\ M_N(\cm)$.\
Because $\ \cm $\  is a II$_1$ factor
all such possibilities arise.
That is, given a projection $\ q$\ in $\ M_N(\cm)$
we can choose $\ N$\  large enough so that $\ $trace$(q_i) \ge
N^{-1}$trace$(q)$,\  for each $\ i.$\  Then there is a
natural partial matrix unit system
$\{f_{k,l}: (k,l) \in E(G)\}$,\ in $\ M_N(\ca)$,\ with
the elementary property above,
such that \ trace$(f_{k,k}) = $\ trace$(q)$\ for all $\ k.$\
If \ trace$(q) = \ga $\ then denote the equivalence class
of these digraph algebras (with $\ H = G$) \ by $\ [A_\ga]$.\

Let $\ f$\  be a partial isometry in $\ M_N(\ca)$,\  for some $\ N$,\
which normalises the subalgebra
$\ M_N(\ca \cap \ca^*)$.\  Then $\ f$\  is elementary in the sense above.
The principle involved here
is that if a partial isometry of the form
\vspace{.3in}

\begin{center}
$
\left[
\begin{array}{cccccc}
0 & 0 & v & 0 & 0 & 0 \\
0 & 0 & 0 & 0 & 0 & w \\
0 & 0 & 0 & 0 & 0 & 0 \\
0 & 0 & 0 & 0 & 0 & 0 \\
0 & 0 & 0 & 0 & 0 & 0 \\
0 & 0 & 0 & 0 & 0 & 0
\end{array}
\right]
$
\end{center}

\noindent normalises the block diagonal algebra of matrices
\vspace{.3in}

\begin{center}
$
\left[
\begin{array}{cccccc}
a & b & 0 & 0 & 0 & 0 \\
c & d & 0 & 0 & 0 & 0 \\
0 & 0 & e & f & 0 & 0 \\
0 & 0 & g & h & 0 & 0 \\
0 & 0 & 0 & 0 & i & j \\
0 & 0 & 0 & 0 & k & l
\end{array}
\right]
$
\end{center}

\noindent then $\ v$\  or $\ w$\  is equal to zero.
It follows that if
$\ A \subseteq M_n(\ca)$\ is an $\ \ca \cap \ca^*$-normalising
digraph algebra for $\ \ca$\
then the partial  matrix unit system for $\ A$\  is
equivalent, by a unitary in $\ M_n(\ca \cap \ca^*)$,\
to a direct sum of subsystems for
the algebras $\ A_\ga$\  identified above.

We now have the identification

\[
H_n(\ca;\cc) = (\sum_{\alpha \in \IR_+} \oplus H_n ([A_\alpha]))/Q_n.
\]

\noindent Identify each group $\ H_n ([A_\alpha])$\  with
$\ H_n(\gD(\cl))$.\
As in the last proof we may replace
$\ Q_n$\ by the subgroup
corresponding to the relations of induced decompositions.
This is the subgroup of
$\ \sum_{\alpha \in \IR_+} \oplus H_n ([A_\alpha])$\
generated by elements of the form

\[
\sum_{\alpha} \oplus  \delta_{\beta, \alpha}g \ - \sum_\alpha \oplus
\delta_{\beta_1,\alpha}g - \sum_\alpha \oplus \delta_{\beta_2,\alpha}g
\]

\noindent where $\ g \in H_n (\Delta(\cl)), \ \beta = \beta_1 + \beta_2$ \  and
where
$\ \delta_{\beta, \alpha}$\  is the Kronecker delta. It follows that
$\ H_n(\ca;\cc) = H_n (\Delta({\cal L}))
\otimes_{\ZZ} \IR $ as desired.\

(ii) The proof of (ii) is very similar to the proof of Theorem 2.1 and so
we omit it. The regularity hypothesis for $\ \cc$\ is necessary because
of the existence of singular masas, that is, masas with trivial normalisers.
\end{mproof}
\vspace{.3in}

\noindent {\bf Example 2.3}
To see that the formula of Theorem 2.2 (i) is not valid
when $\ \cm$\  is
the $\ I_\infty$\   factor let $\ \ca = A(D_4))$\
be the subalgebra of $\ M_4 (\IC)$\  spanned by the matrix units $\ e_{13},
e_{14},  e_{23},  e_{24}$\  \ and the standard diagonal subalgebra $\ \IC^4$.\
This is
the standard example of a matrix algebra with nontrivial Hochschild
cohomology and the last theorem shows that $\ H_1 (\ca ; \cc) = \ZZ$.\
Let $\ \cal B$\  be the operator algebra on $\ \IC \oplus (\IC^4 \otimes {\cal
H})$,\  where $\ \cal H$\  is an infinite dimensional Hilbert space,
consisting of operators of the form

\[
{\ \left( \begin{array}{cc} \lambda &
* \\ 0 & a \end{array} \right),\ }
\ \  \mbox{with} \ \lambda \in \ \IC \ \ \mbox{and}\  \  \ a \in \ca \otimes
\cl({\cal H}).\
\]

\noindent In the terminology of Gilfeather and Smith
\cite{gil-smi-1} this is the {\em cone algebra} of
$\ \ca \otimes \cl({\cal H})$.\  The algebra $\ \cb$\  is the reflexive
operator algebra determined by a finite commutative projection lattice,
with five atoms $\ q_1, q_2, q_3, q_4, q_5$\  whose associated digraph
for the order
$\ \ll$,\  is a 4-cycle (for $\ q_1, q_2, q_3, q_4$\ ) with an added vertex
(for $\ q_5$\ )
which receives four directed edges from each of the vertices of the
4-cycle.
Here $\ q_5$\  is the rank one projection onto the one dimensional
summand, and $\ q_i$\  is $\ e_{i,i} \otimes I_{\ch}$,\  for $\ 1 \le i \le
4$.\
Although $\ H_1(\Delta(\cl))$\  is zero, clearly, the basic 1-cycle
$\ (e_{1,1} \otimes I_{\ch}, e_{2,2} \otimes I_{\ch},
e_{3,3} \otimes I_{\ch}, e_{4,4} \otimes I_{\ch})$\  gives a generator for
$\ H_1(\cb ; \cb \cap \cb^*)$.\
\vspace{.3in}

The next identification is in the context of limit algebras, one of
our key motivating contexts for the formulation of stable homology.

Let
$\ \displaystyle{\ca = \lim_\to (A(G_k),\phi_k)}$\
be the Banach algebra direct limit of a direct system
of digraph algebras $\ A(G_k)$\  with star-extendible injections
$\ \phi_k : A(G_k) \to A(G_{k+1})$\  which map standard matrix units
to sums of standard matrix units. In particular,
such maps are {\em regular} in the sense of the next section.
For each $\ n \ge 0$\
 there is a natural induced group homomorphism
\[
(\phi_k)_* : H_n(\Delta(G_k)) \to H_n(\Delta(G_{k+1}))
\]
and an associated
direct limit abelian
group
\[\displaystyle{\lim_\to (H_n(\Delta(G_k)), (\phi_k)_*)}.
\]

\noindent Such limit groups have appeared in \cite{dav-scp}
and \cite{scp-book}.
Let
$\ \displaystyle{\cd = \lim_\to (\ \IC^{|G_k|},\phi_k)}$\
be the abelian \cssa \ of $\ \ca$,\  where \ \ $\ \IC^{|G_k|}$\ \
is the standard diagonal subalgebra of $\ A(G_k)$.\
\vspace{.3in}

The following matricial variant of a fundamental fact for normalising
partial isometries in AF \csa s will be needed.
The scalar case appears as Lemma 5.5 of \cite{scp-book}.
\vspace{.3in}

\noindent {\bf Lemma 2.4}\ \ \it Let $\ \cb, \cd$\ be as above and let $\ f $\
be a partial isometry
in $\ \cb \ot M_m$\  which normalises $\ \cd \ot M_m.$\ Then
$\ f = dw$\ where $\ d$\ is a partial isometry in $\ \cd \ot M_m$\ and
$\ w$\ is a partial isometry in $\ \cb_k$,\ for some $\ k,$\ which normalises
the diagonal subalgebra $\ \cd_k.$
\rm
\vspace{.3in}

\begin{mproof} Let $\ \tilde{\cb_k}$\ be the algebra generated by
$\ \cb_k$\ and $\ \cd$\ and let
$\ P_n: \cb \to \tilde{\cb_k}$  be the
natural projections, as given in Chapter 4 of \cite{scp-book}
for example. In particular $\ P_n$\ is the pointwise limit of maps
$\ P_{n,r}, \ \ r = 1, 2, \dots$,\  each of which has the form
$\ P_{n,r}(b) = p_1bp_1 + \dots + p_rbp_r   $\ for some
family of orthogonal projections in $\ \cd$. This
property shows that if $\ v \in \cb$\ is a partial isometry normalising
$\ \cd$\ then so too is each operator $\ P_{n,r}(v)$,\  and hence so too is
$\ P_n(v)$\  itself. It
follows that the map $\ P_n \ot $Id\ $: \cb \ot M_m \to
\tilde{\cb_n} \ot M_m$\ is defined in such a way that it follows that
$\ (P_n \ot $Id$)(f)$\ is also a
partial isometry which normalises$\ \cd \ot M_m$.\

We can now argue exactly as in the proof in \cite{scp-book}
for the scalar case $\ n = 1.$

Let $\ q_n$\ be the range projection of
$\ (P_n \ot $Id$)(f).$\  Then \ $q_n$\ is a Cauchy sequence
of projections in $\ \cd \ot M_m$\  converging to $\ ff^*.$ \
Since $\ \cd$\  is abelian
it follows that there exists $\ n_0$\ such that $\ q_n
\in \tilde{D_{n_0}   } \ot M_m$\  for all $\ n$.\
The lemma is straightforward in the special case $\ \cb = \tilde{\cb_t}$,
\ and so it will be sufficient to prove that
$\ (P_n \ot $Id$)(f) \in
\tilde{\cb_{n_0}} \ot M_m
\ $ for all $\ n$,\ since from this it follows
that $\ f \in \tilde{\cb_{n_0}} \ot M_m$.

Write $\ (P_n \ot $Id$)(f) = (P_{n_0} \ot $Id$)(f)  + z$.\  Then \
$z = \sum c_i e_i$,\  a finite sum with
coefficients $\ c_i $\ in $\ \cd \ot M_m$\ and where each $\ e_i$\ is a
standard matrix unit for $\ \cb_n$\ which is not
subordinate to a standard matrix unit for $\ \cb_{n_0}$.\
It follows that
\ $(P_n \ot $Id$)(z(P_{n_\circ} \ot $Id$) (f)^*) = 0$\ for $\ n \ge n_0.$\
Thus
\[
ff^* = ((P_{n_\circ}  \ot \mbox{Id})(f) + z) ((P_{n_\circ} \ot \mbox{Id}) (f) +
z)^*
\]
\[
= q_{n_0} + zz^* + z((P_{n_\circ} \ot \mbox{Id})(f))^* + (P_{n_\circ} \ot
\mbox{Id})(f)z^*
\]
and so
\[
(P_n \ot \mbox{Id})(ff^*) = (P_n \ot \mbox{Id})(ff^*) + (P_n \ot
\mbox{Id})(zz^*).
\]
Thus $\ (P_n \ot $Id$)(zz^*) = 0 $.\
\ Let $\ n \to \infty$\  and we obtain \ $(P \ot $Id$)(zz^*) = 0$. Since
\ $(P \ot $Id)\  is a faithful expectation \ $\ z = 0$\  as desired.
\end{mproof}
\vspace{.3in}

\noindent {\bf Theorem 2.5}\ \  \it Let $\ \ca$\  be the operator algebra
$\ \displaystyle{ \lim_\to (A(G_k),\phi_k)}$\
with regular embeddings and diagonal
subalgebra $\ \cd$,\  as above. Then, for each $\ n \ge 0,$\
the stable homology
group $\ H_n(\ca;\cd)$\  is isomorphic to the limit homology group
\ \ $\ \displaystyle{\lim_\to (H_n(\Delta(G_k))), (\phi_k)_*)}.$\
\rm

\begin{mproof} Let $\ A \subseteq M_\infty(\ca)$\  be a $\ \cd$-normalising
 digraph algebra for  $\ \ca$\  with a partial matrix unit system
$\ \{f_{i,j} : (i,j) \in I_A\}$\  which generates a full matrix unit system
$\ \{f_{i,j} : 1 \le i, j \le m \}$\  in $\ M_\infty(\cb)$,\
where $\ \cb$\  is the AF \csa \ generated by $\ \ca$.\  Without loss of
generality
assume that the digraph of $\ A$\  is connected. From
Lemma 2.4 it follows that
the full system
$\ \{f_{i,j}\}$\  is unitarily equivalent, by a unitary in $\ M_\infty(\cd)$\
to a system $\ \{g_{i,j}\}$\  where, for some integer $\ k > 0$,\  each
$\ g_{i,j}$\   is a sum of the standard matrix units of the subalgebra
$\ M_k(C^*(A(G_k)))$\  of $\ M_\infty(\cb)$.\  Here we identify $\ A(G_k)$\
and
 its generated \csa \ with its image in$\ \ca$\  and $\ C^*(\ca)$\
respectively.
It follows, as in the proof of Theorem 2.1, that
\[
H_n(\ca;\cd) = (\sum_k \oplus H_n([M_k(A(G_k))]))/Q_n
\]
and so
\[
H_n(\ca;\cd) =  (\sum_k \oplus H_n([A(G_k)]))/Q_n.
\]
Furthermore in the second quotient expression we may assume that $\ Q_n$\
is the set of relations for  the standard inclusions and induced
decompositions amongst the set of digraph algebras $\ M_k(A(G_k))$.\

Let $\ \eta$\  be the natural group homomorphism from the direct limit
group $\ \cg$\  say, to $\ H_n(\ca;\cd)$.\
This is well-defined, because the relations $\ Q_n$\  include those
relations coming from the given injections $\ \phi_k$.\
On the other hand,
suppose that $\ h \in H_n([A(G_k)])$\  and $\ h \in Q_n$.\
Then there exists $\ k_1 > k$\  so that $\ g$\  is a finite sum of terms
of the form\ \
$\ g - \theta(g)$\  \ \ and \ \
$\ g - \theta^\prime(g) - \theta^{\prime \prime}(g)$\ \ \
associated with the given inclusions
$\ A(G_p) \to A(G_{k_1}),$\  for $\ 1 \le p \le k$.\
Thus, viewed as a member of $\ A(G_k)$,\  $\ g$\  is the zero element.
Hence $\ \eta$\  is injective and surjective.
\end{mproof}
\vspace{.3in}

The following examples
can be obtained readily with the help
of the theorems above.
\vspace{.3in}

\noindent {\bf Example 2.6} \  Let $\ \ca $\ be a strongly
maximal triangular subalgebra of the AF \csa \ $\ \cb$.\
(See \cite{ppw1}, \cite{scp-book}.)
Then $\ H_1(A(D_4) \ot \ca) = K_0(\ca \cap \ca^*).$\
Here the unique masa $\ \ca \cap \ca^*$\ is understood and suppressed from the
notation.

On the other hand
let $\ A(D_4) \star \ca$\ be the lexicographic product (cf. \cite{scp-lex-1})
given by

\[
(A(D_4) \cap A(D_4)^*)
 \ot \ca \ +\  A(D_4)^0 \ot \cb,
\]

\noindent where $\ A(D_4)^0$\ is the kernel of the diagonal expectation
onto the diagonal algebra \ $ A(D_4)\cap A(D_4)^*$.
This algebra is triangular, with a unique masa, and
 $\ H_1(A(D_4) \star \ca) = K_0(\cb).$\
\vspace{.3in}

\noindent {\bf Example 2.7} \ Let

\[
\phi_k : A(D_4) \ot (M_{3^k} \op M_{3^k}) \to
A(D_4) \ot (M_{3^{k+1}} \op M_{3^{k+1}})
\]

\noindent be the embedding $\ \phi \ot {id}_{M_{3^{k-1}}}$,\ where
$\ \phi$\ is the embedding given before Definition 3.3.
(Identify $\ (M_{3^k} \op M_{3^k})$ \
with $\ ((M_{3} \op M_{3}) \ot M_{3^{k-1}}$ \
etc.) Let $\ \ca$\ be the associated unital digraph limit algebra.
Then, with respect to the natural diagonal subalgebra $\ \cc$,\
\[
H_1(\ca;\cc) \ = \
\displaystyle{\lim_\to (\ZZ^2, \small \left[ \begin{array}{cc}
2 & 1\\1 & 0
\end{array} \right]\rm )}
\ =\  \ZZ^2.
\]
\vspace{.3in}

There are a number of interesting connections between
$\ H_n(\ca;\cc)$\ and automorphism, derivations and Hochschild
cohomology.
The following theorem is an example of this.
Related assertions (with similar proofs)
can be found in  \cite{dav-scp} \cite{scp-simplicial} and \cite{scp-spectral}.

Let $\ \cc$\ be a maximal abelian algebra in $\ \ca$\  and write
\ Aut$_{\cc}(\ca)$\  for the corresponding group of Schur automorphisms
of $\ \ca.$\  This is the group of  automorphisms $\ \ga$\ for
which $\ \ga(c_1ac_2) = c_1\ga(a)c_2$\ \ for all $\ c_1, c_2$\ in $\ \cc$\
and $\ a$\ in $\ \ca$.\  If $\ A$\ is a $\ \cc$-normalising digraph algebra
for $\ \ca$\
then write $\ {\tilde A} $\ for the algebra generated
by $\ A$\ and $\ \cc$.\  We say that a Schur automorphism
is {\it locally} \ {\it $\ \cc$-inner}\ if the restriction to each such
subalgebra
$\ {\tilde A} $\  is inner.
\vspace{.3in}

\noindent {\bf Theorem 2.8}\ \ \it Let $\ \cc$\ be a maximal abelian subalgebra
of the operator algebra $\ \ca$\ \ and suppose that $\ H_1(\ca;\cc) = 0.$\ Then
every
Schur automorphism in \ Aut$\ _{\cc}(\ca)$\ is locally \  $\ \cc$-inner.\rm
\vspace{.3in}

A locally \  $\ \cc$-inner automorphism need not be inner
even for approximatelty finite C*-algebras
and their regular subalgebras. (See, for example, Remark 2 of
\cite{scp-outer}.) Nevertheless such automorphisms are often approximately
inner
in the sense of being approximable in the point-norm topology
by inner (Schur) automorphisms. Thus, in rough parallel with the weakly closed
theory developed in \cite{scp-spectral}, it seems to be the case that
there is a close connection between stable homology with respect to regular
maximal abelian self-adjoint subalgebras, and the (norm) {\it essential}
Hochschild
cohomology arising when boundaries are replaced by their point norm closures
- the essential boundaries.

\newpage

\section {\bf Regular inclusions and $\ K_0 \oplus H_*$\  uniqueness}

The following distinguished class of embeddings is studied in
\cite{scp-book}, \cite{scp-simplicial} and \cite{scp-k0}.
\vspace{.3in}

\noindent {\bf Definition 3.1}\ \cite{scp-k0} \ A
star-extendible algebra homomorphism between digraph algebras is said
to be {\em regular} if it is (inner) unitarily equivalent to a direct
sum of multiplicity one star-extendible embeddings.
\vspace{.3in}

A multiplicity one star-extendible embedding $\ A(G) \to A(H)$\  is
a restriction of a star homomorphism $\ C^*(A(G)) \to
C^*(A(H))$\  which is of multiplicity one. In particular every
star homomorphism between self-adjoint digraph algebras is
automatically regular.
On the other hand there are, in general, a myriad of star-extendible
homomorphisms between digraph algebras, and the regular
embeddings form the most natural subclass. Between two digraph
algebras there are only finitely many  (inner)
unitary equivalence classes of  regular homomorphisms, and, for
elementary algebras, these classes may be represented by diagrams
at the level of digraphs. Bratteli diagrams form a  degenerate case.
The terminology
"regular" is used because direct systems of  regular
embeddings provide limit algebras possessing a distinguished
maximal abelian self-adjoint subalgebra
which is regular in the usual sense that the normaliser of the masa
generates the algebra.

An important aspect of regular morphisms is that they are the correct
class of maps to  consider with regard to the functoriality of stable
homology;
each regular homomorphism \ $\phi : A(G) \to A(H)$ \ induces a group
homomorphism \ $\phi_* : H_n(A(G)) \to H_n(A(H))$. \ Here we have
written
$\ H_n(A(G))$\ for $\ H_n(A(G);C)$\
where $\ C$\  is any maximal abelian subalgebra
of $\ A(G)$.\  This is a well-defined move since each such masa
is unique up to inner unitary equivalence.

If we focus attention on a stable family of digraph
algebras of the form $\ A(G) \otimes M_n$,\  \ $\ n = 1,2....,$\  \ where \ $\
G$\  \ is a
fixed digraph, then the following class of regular embeddings is
particularly natural. As we shall see these
rigid embeddings appear naturally in the construction of limit algebras with
interesting homology. Furthermore,
for various stable families we can
classify associated rigid inclusions in terms of the induced map
on $\ K_0 \oplus H_*$.\
\vspace{.3in}

\noindent {\bf Definition 3.2.} \ \cite{scp-book}\ (i) Let $\ G$\  be a
connected
digraph. A {\em rigid embedding}
$\ A(G) \otimes M_n \to A(G) \otimes M_m$\  \ is a regular embedding which
is unitarily equivalent
to a direct sum of embeddings
$\ \theta \ot \psi $\
where $\ \psi : M_n \to
M_m$\  is a multiplicity one \csa \ algebra injection
and $\ \theta: A(G) \to A(G)$\  is an
automorphism induced by a digraph automorphism.

(ii) A general rigid embedding $\ A(G) \otimes B_1 \to A(G) \otimes B_2$,\
with $\ B_1, B_2$\  finite-dimensional $\ C^*$\ -algebras, is a
star-extendible embedding for which the
partial embeddings are rigid.
\vspace{.3in}

The unitary equivalence class of a rigid embedding can be indicated
by a (unique) labelled Bratteli diagram in which each edge from a vertex $\ i$\
of level one to vertex $\ j$\  of level two indicates a multiplicity one
partial rigid embedding, and the labelling of the
edge indicates the particular automorphism $\ \theta$\  used in the embedding.

For example, let $\ \theta_1 $\  and $\ \theta_3$\  be the identity
and rotation automorphisms of $\ A(D_4)$,\
and let $\ \theta_2 $\  and $\ \theta_4$\   be the two reflections. The diagram

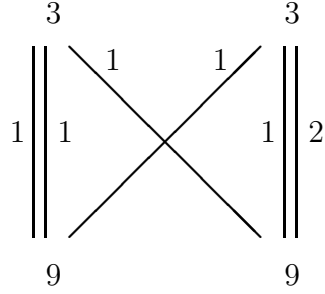
\begin{center}
\setlength{\unitlength}{0.0125in}%
\begin{picture}(125,119)(75,700)
\thicklines
\put(100,800){\line( 1,-1){ 80}}
\put(100,720){\line( 1, 1){ 80}}
\put(190,800){\line( 0,-1){ 80}}
\put( 85,800){\line( 0,-1){ 80}}
\put(195,800){\line( 0,-1){ 80}}
\put(190,810){\makebox(0,0)[lb]{\raisebox{0pt}[0pt][0pt]{\twlrm 3}}}
\put( 95,760){\makebox(0,0)[lb]{\raisebox{0pt}[0pt][0pt]{\twlrm 1}}}
\put( 90,800){\line( 0,-1){ 80}}
\put(180,760){\makebox(0,0)[lb]{\raisebox{0pt}[0pt][0pt]{\twlrm 1}}}
\put(190,700){\makebox(0,0)[lb]{\raisebox{0pt}[0pt][0pt]{\twlrm 9}}}
\put(115,790){\makebox(0,0)[lb]{\raisebox{0pt}[0pt][0pt]{\twlrm 1}}}
\put(160,790){\makebox(0,0)[lb]{\raisebox{0pt}[0pt][0pt]{\twlrm 1}}}
\put( 75,760){\makebox(0,0)[lb]{\raisebox{0pt}[0pt][0pt]{\twlrm 1}}}
\put(200,760){\makebox(0,0)[lb]{\raisebox{0pt}[0pt][0pt]{\twlrm 2}}}
\put( 90,810){\makebox(0,0)[lb]{\raisebox{0pt}[0pt][0pt]{\twlrm 3}}}
\put( 90,700){\makebox(0,0)[lb]{\raisebox{0pt}[0pt][0pt]{\twlrm 9}}}
\end{picture}
\end{center}

\noindent indicates the rigid embedding
\[
\phi : A(D_4) \otimes (M_3\ \oplus \ M_3)
\to A(D_4) \otimes (M_9\ \oplus \ M_9).
\]
where all multiplicity one component embeddings have the identity
automorphism excepting that for the edge labelled with a 2,
which is the reflection $\ \theta_2$.
One can verify that (with natural identifications of the homology groups)
$\ \phi$\  induces maps
$\ H_0\phi : \ZZ^2 \to  \ZZ^2$\  and $\ H_1\phi :  \ZZ^2 \to  \ZZ^2$\
given by
\[
H_0\phi =
\left[ \begin{array}{cc}
2 & 1\\1 & 2
\end{array} \right] \ \ \hbox{and} \ \
H_1\phi =
\left[ \begin{array}{cc}
2 & 1\\1 & 0
\end{array} \right].
\]
\vspace{.3in}

\noindent {\bf Definition 3.3}\  \ (i) A {\it cycle algebra}, or {\it
2m-cycle digraph algebra},
is a digraph algebra of the form $\ A(D_{2m}) \otimes B$\
where $\ D_{2m}$\ is the $2m$-cycle digraph and where
$\ B$\  a
finite-dimensional $\ C^*$\ -algebra.\\
(ii) If $\ A_1 \subseteq A_2$\  are $\ 2m$-cycle digraph algebras, then
the inclusion is said to be
{\em rigid} if the inclusion map is a rigid embedding.
\vspace{.3in}

The following proposition is elementary but it is a direct
counterpart to the important  fact
that inclusions of finite-dimensional $\ C^*$-algebras are determined
up to inner conjugacy
by their induced $\ K_0$\  maps.
\vspace{.3in}

\noindent {\bf Proposition 3.4}\it \  \ A rigid embedding between cycle
algebras is
determined up to inner
unitary equivalence by the induced maps between the $\ K_0$\  groups
and the first stable homology groups.\rm

\begin{mproof} Let $\ D_{2m}$\  be a $\ 2m$-cycle digraph
with receiving vertices
labelled $\ v_1, v_3,\dots,v_{2m-1}$\  and emmitting vertices
$\ v_2, v_4, \dots ,v_{2m}$.\  Let $\ \theta_1, \theta_3, \dots
,\theta_{2m-1}$\
be the rotation automorphisms of $\ D_{2m}$\  such that $\ \theta_j(v_1)
=v_j$,\  and let $\ \theta_2, \theta_4, \dots,  \theta_{2m}$\
be the reflection automorphisms $\ \theta_{2j} = \eta \circ \theta_{2j-1},
1 \le j \le m$,\  where $\ \eta$\  is the reflection fixing $\ v_1$.\
Write $\ \theta_k$\  also for the automorphisms of $\ A(D_{2m})$\
induced by these graph automorphisms.

A rigid embedding $\ \phi : A(D_{2m}) \otimes M_p \to A(D_{2m})
\otimes M_q$\
is unitarily equivalent to the direct sum
$\ r_1\theta_1 + \dots + r_{2m}\theta_{2m}$\  where we abuse
notation and write $\ r_k\theta_k$\  for the orthogonal direct sum
of $\ r_k$\  copies of the embeddings $\ \theta_k \otimes id$.\  Clearly the
$\ 2m$-tuple $\ {r_1, \dots ,r_{2m}}$\  is a complete invariant for
the unitary equivalence class of $\ \phi$.\
It will be enough to show that the inner equivalence
class of $\ \phi$\ is determined by this  2m-tuple.

The map $\ K_0\phi,$\  under the natural identification of the $\ K_0$\
groups,
has the form $\ X + JY$\  where $\ X = X(r_1, r_3,\dots,r_{2m-1})$\  is the
Laurent matrix

\[
\left [\begin {array}{ccccccccc}
r_1&0&{  r_{2m-1}}&.&.&.&.&r_3&0
\\0&r_1&0&.&.&.&.&0&r_3
\\r_3&0&r_1&.&.&.&.&r_5&0
\\0&r_3&0&.&.&.&.&0&r_5
\\.&.&.&.&.&.&.&.&.
\\r_{2m-1}&0&.&.&.&.&.&r_1&0
\\0&r_{2m-1}&.&.&.&.&.&0&r_1
\end {array}\right ],
\]

\noindent where $\ Y$\  is the Laurent matrix $\ X(r_2,r_4,\dots,r_{2m})$,\
and where $\ J$\  is the matrix

\[
\left [\begin {array}{cccccc}
1&0&.&.&0&0
\\0&0&.&.&0&1
\\0&0&.&.&1&0
\\.&.&.&.&.&.
\\0&0&.&.&.&.
\\0&1&.&.&0&0
\end {array}\right ].
\]

\noindent On the other hand the map $\ H_1\phi : \ZZ \to \ZZ$,\  under the
natural identification of the $\ H_1$\  groups, is $\ H_1\phi = [\delta]$\
where
\[
\delta = (r_1 + r_3 + \dots + r_{2m-1}) -(r_2 + r_4 + \dots + r_{2m}).
\]
The proposition will be proven if we show that the two matrices $\ K_0\phi$\
and $\ H_1\phi$\
determine the coefficients $\ r_1, \dots ,r_{2m}$.\
To this end
let $\ \pi : M_{2m} \to M_{2m}$\  be the natural projection
onto the Laurent matrices obtained by averaging
the $\ 2m$\  entries of each of the $\ m$\  odd "diagonals"
and replacing the other diagonals with zeros.
Note that if $\ X$\  is a Laurent matrix then $\ \pi(JX)$\  is a multiple of
the "all ones" matrix $\ Z = X(1,1,\dots,1)$.\
It follows that application of $\ \pi$\  to the matrix $\ X + JY$\
determines the components $\ X, Y$\  up to a multiple of $\ Z$.\
That is,  the ordered sets
$\ \{r_1,r_3,\dots,r_{2m-1}\}, $\ \ and$\ \ \{r_2,r_4,\dots,r_{2m}\}$\
are determined up to a common additive constant. But
now the fact that the difference $\ \delta$\
is given by $\ H_1\phi$\  leads to the
determination of $\ r_1,\dots,r_{2m}$.\
\end{mproof}
\vspace{.3in}

\noindent {\bf Corollary 3.5} \ \it Let $\ A_1, \ A_2, \ and
\ A$\  be $\ 2m$-cycle digraph algebras
with  $\ A_1 \subseteq A, \ A_2 \subseteq A$\  where the inclusions are rigid.
Then $\ A_1$\  and $\ A_2$\  are inner conjugate if and only if the inclusion
maps
induce the same maps between the $\ K_0$\  groups and between
the first stable homology groups. \rm
\vspace{.3in}

\noindent {\bf Definition 3.6}\ \ Let $\ G$\  be a digraph and let $\ \Theta$\
be a
subset of
\ Aut$(G)$.\
Then $\ \Theta$\  is said to have the
{\it $\ K_0 \op H_1-$}{\it uniqueness property} \ if the rigid embeddings from
$\ A(G) \otimes M_p$\   to $\ A(G) \otimes M_q$\  which are associated with
$\ \Theta$\  are determined up to inner conjugacy by the induced maps on
$\ K_0$\  and $\ H_1$.\  The $\ K_0 \op H_*$\  -{\em uniqueness property} is
defined similarly.
\vspace{.3in}

As part of the general homology programme for limit algebras indicated in
\cite{scp-book} it is of interest to determine contexts (G,$\ \Theta$\ )
which have the $\ K_0 \op H_*$\  - uniqueness property. This
gives a starting point for classifications of \nsa \ limit
algebras in the style of Elliott's classification of AF \csa s.
\vspace{.3in}

\noindent {\bf Example  3.7 Suspensions}
\ \
Let \ $K_n^i,\  i  = 1, 2,$ \ be complete digraphs on n vertices.
Define the {\em n-point suspension} of the digraph algebra $\ A = A(G)\ $
to be
the digraph algebra  $\ S_nA\ $ with graph $\ S_nG\ $
where the vertex and edge sets
are given by

\[
V(S_nG) = V(K_n^1)\  \cup V(K_n^2) \ \cup V(G),
\]

\[
E(S_nG) = E(K_n^1) \ \cup E(K_n^2) \ \cup E(G) \cup E
\]

\noindent where $\ E =
\{(v^i,w) : w \in V(G), v^i \in V(K_n^i), i = 1,2\}
$.
Let $\ G_1, G_2$\ \ be connected.
A regular
embedding $\phi : A(G_1) \to A(G_2)$ of multiplicity $r$
induces a natural regular embedding $S_k\phi : S_k(A(G_1)) \to S_{kr}(A(G_2))$
which respects the north pole and south pole summands of the
suspension algebras. This suspended embedding
is uniquely determined up to inner conjugacy.
{}From simplicial homology theory it follows that for each order $t$
the suspended embedding $\ S_k\phi \ $
induces a homomorphism of the stable homology groups of order $t+1$,
and this homomorphism may be identified with the homomorphism of the
homology groups of order $t$ induced by $\ \phi. $
It follows that the homological classifications in this paper
of various families of embeddings admit immediate
higher order extensions to the classification
of the associated pole preserving embeddings of the suspension algebras.
\vspace{.3in}

\noindent {\bf Example 3.8 Discrete Tori.}\ The
discrete tori algebras are the digraph
algebras
\[
A(D_{2m_1}) \otimes \dots \otimes A(D_{2m_s})
\]
whose underlying digraphs are the direct products of cycle digraphs.
The full group of rigid automorphisms of these algebras fails to have
the $\ K_0 \op H_*$- uniqueness property.
To see this consider the
rigid embeddings

\[
\phi,\ \psi :\ A(D_4) \otimes A(D_4) \to \ A(D_4) \otimes A(D_4)
\otimes M_{12}
\]

\noindent  given by

\[
\phi = ((2\theta_1 \oplus \theta_3) \otimes (\theta_1 \oplus \theta_3))
\oplus ((\theta_1 \oplus 2\theta_3) \otimes (\theta_2 \oplus \theta_4)),
\]
\[
\psi = ((\theta_1 \oplus 2\theta_3) \otimes (\theta_1 \oplus \theta_3))
\oplus ((2\theta_1 \oplus \theta_3) \otimes (\theta_2 \oplus \theta_4)).
\]

\noindent Then $\ K_0\phi$\  and $\ K_0\psi$\  coincide with
$\ 3X \otimes X$,\ \  where $\ X$\  is the "all ones" matrix
$\ X(1,1,1,1).$\  Also one can verify that
$\ H_0\phi = H_0\psi = [12],\
H_2\phi = H_2\psi = [0], $\ \ and $\  H_1\phi = H_1\psi = 0,$\   the zero map
from
$\ \ZZ^2$\  to $\ \ZZ^2$.\  Thus $\ (K_0 \oplus H_*)\phi =
(K_0 \oplus H_*)\psi$\ \ \ and yet
the injections are not inner conjugate.
\vspace{.3in}

\noindent {\bf Example 3.9 \ \ The Cube Algebra.}\ \
Define the {\it cube algebra} to be the digraph algebra
 in $\ M_8$\  which is associated
with the following digraph, which we denote as  $\ \ Cu$.

% cube1.tex
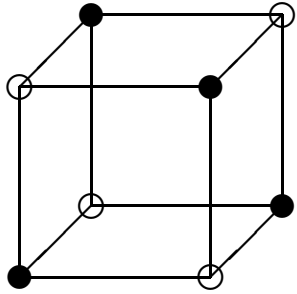
\begin{center}
\setlength{\unitlength}{0.0125in}%
\begin{picture}(120,120)(125,685)
\thicklines
\put(130,690){\circle*{10}}
\put(160,800){\circle*{10}}
\put(210,770){\circle*{10}}
\put(210,690){\circle{10}}
\put(240,800){\circle{10}}
\put(130,770){\circle{10}}
\put(240,720){\circle*{10}}
\put(160,720){\circle{10}}
\put(160,720){\line(-1,-1){ 30}}
\put(160,800){\line( 1, 0){ 80}}
\put(240,800){\line( 0,-1){ 80}}
\put(240,720){\line(-1, 0){ 80}}
\put(160,720){\line( 0, 1){ 80}}
\put(160,800){\line(-1,-1){ 30}}
\put(130,770){\line( 1, 0){ 80}}
\put(210,770){\line( 0,-1){ 80}}
\put(210,690){\line(-1, 0){ 80}}
\put(130,690){\line( 0, 1){ 80}}
\put(240,800){\line(-1,-1){ 30}}
\put(240,720){\line(-1,-1){ 30}}
\end{picture}
\end{center}

\noindent This may be regarded as a three dimensional variant of the 4-cycle
graph which appears on each face of the cube.
The full automorphism group Aut$\ (Cu)$\  has 24 elements corresponding
to the 24 permutations
of the receiving vertices. Note that there is a unique
directed graph automorphism of $\ Cu$\  for each such permutation.
Thus Aut$\ (Cu)$\  has order 24, and a general rigid embedding
$\ \phi : A(Cu) \otimes M_n \to A(Cu) \otimes M_m$\  has an inner unitary
conjugacy class which is determined by the ordered set $\
\{r_1,\dots,r_{24}\}$\
corresponding to the multiplicities of the types of partial rigid embeddings.
Furthermore,
it follows that in the direct sum decomposition
\[
K_0\phi = K_0^r\phi \oplus K_0^e\phi,
\]
corresponding to the receiving and emmitting summands, the linear system in
the unknowns  $\ \{r_1,\dots,r_{24}\}$\  arising from the equation
$\ K_0^r\phi = K_0^r\psi,$\  with $\ \psi$\  given, has the same rank as the
system
for the full equation $\ K_0\phi = K_0\psi$.\
Thus, knowledge of the $\ 4$\  by $4$\
matrix $\ K_0\phi$\  leads to 16 equations for
$\ \{r_1,\dots,r_{24}\}$.\  We have $\ H_1(A(Cu) \otimes M_n) =
\ZZ^5$,\  and so 25 more equations are provided
by $\ H_1\phi$\  giving a  system of  41 linear equations in 24 unknowns.
Curiously, (computer assisted) calculation shows
that the coefficient matrix of this system has rank 23 and so the
full automorphism group for the cube algebra just misses having the
$\ K_0 \op H_*$\  - uniqueness property. This can be seen directly
by considering the multiplicity 12 embedding which is a direct sum of
the rotations and the multiplicity 12 embedding which is the direct
sum of the rest. Both induce the zero map on $\ H_1$\  and both have the
same $\ K_0$\ map.

On the other hand,  proper subgroups of
Aut$\ (Cu)$\  do have this uniqueness
property. In particular, this is the case for the
group of 12 orientation preserving
symmetries of the cube digraph.
Calculation shows that the coefficient matrix in this simpler case is
the following.
\vspace{.3in}

\begin{center}
\bf Coefficient Matrix arising from
Rotations of $\ Cu$\ \rm
\tiny

\[
\left [\begin {array}{cccccccccccc} 1&0&0&0&0&1&0&0&0&0&1&0\\0&0&0&1&0
&0&0&0&1&0&0&1\\0&1&0&0&1&0&0&0&0&1&0&0\\0&0&1&0&0&0&1&1&0&0&0&0\\0&0&0
&1&0&0&1&0&0&1&0&0\\0&1&0&0&0&0&0&1&0&0&1&0\\1&0&1&0&0&0&0&0&0&0&0&1\\0
&0&0&0&1&1&0&0&1&0&0&0\\0&0&0&0&1&0&0&1&0&0&0&1\\0&0&1&0&0&1&0&0&0&1&0
&0\\0&0&0&0&0&0&1&0&1&0&1&0\\1&1&0&1&0&0&0&0&0&0&0&0\\0&1&1&0&0&0&0&0&
1&0&0&0\\1&0&0&0&1&0&1&0&0&0&0&0\\0&0&0&1&0&1&0&1&0&0&0&0\\0&0&0&0&0&0
&0&0&0&1&1&1\\0&0&-1&1&-1&0&0&0&0&0&1&0\\-1&0&0&0&0&0&0&1&1&-1&0&0\\0&
-1&0&0&0&1&1&0&0&0&0&-1\\1&0&0&0&0&0&0&-1&-1&1&0&0\\0&1&0&0&0&-1&-1&0&0
&0&0&1\\0&1&-1&0&-1&0&1&0&0&0&0&0\\-1&0&0&0&1&0&0&0&0&-1&1&0\\0&-1&0&0
&0&0&0&0&1&1&0&-1\\0&0&1&1&0&0&0&-1&-1&0&0&0\\1&0&0&0&0&-1&-1&1&0&0&0&0
\\1&0&-1&0&-1&0&0&0&1&0&0&0\\-1&0&0&0&0&0&1&0&0&-1&0&1\\0&-1&1&0&0&0&0
&0&0&0&1&-1\\0&1&0&0&0&1&0&-1&-1&0&0&0\\0&0&0&1&1&-1&-1&0&0&0&0&0\\0&0
&-1&0&-1&1&0&0&0&0&0&1\\-1&0&1&1&0&0&0&0&0&-1&0&0\\1&-1&0&0&0&0&0&1&0&0
&0&-1\\0&0&0&0&1&0&0&-1&-1&0&1&0\\0&0&0&0&0&-1&-1&0&1&1&0&0\\0&0&-1&0&
-1&0&0&1&0&1&0&0\\-1&1&0&0&0&1&0&0&0&-1&0&0\\0&-1&0&1&1&0&0&0&0&0&0&-1
\\0&0&0&0&0&0&1&-1&-1&0&0&1\\0&0&1&0&0&-1&-1&0&0&0&1&0\end {array}
\right ]
\]
\end{center}
\rm
\vspace{.3in}

\noindent The rank of the matrix is $\ 12$.\
The submatrix arising from the $\ K_0$\  data
alone is the $\ 16$\  by $\ 12$\  submatrix formed
by the first $\ 16$\  rows, and this has rank $\ 10$.\  Thus,
as in the case of the  cycle algebras,
the stable homology information is really needed.
We have
\vspace{.3in}

\noindent \bf Theorem 3.10 \it \ \ Let $\ \ca $\ be the cube algebra
$\ A(Cu) \ot M_n.$\  Let $\ \cf  $\ be the family of subalgebras of $\ \ca$\
which
are completely isometrically
isomorphic to a cube algebra $\ A(Cu) \ot M_r,$\  for some $\  r,$\
and for which the inclusion map is a rigid embedding associated with
rotations. Then the algebras in $\ \cf  $\  are classified up
to inner conjugacy by the following two invariants.

(i) the inclusion induced map between the scaled $\ K_0$\  groups,

(ii) the induced map between the first stable homology groups.
\rm
\vspace{.3in}

As we can see, even for simple digraph algebras the $\ K_0 \op H_*$\  data can
generate a large system for the
unknown multiplicities
of the components. It is of interest therefore to
discover general combinatorial principles that can assist with rank
determination.

\newpage

\section{Regular Inclusions in AF algebras}

\noindent We now consider regular inclusions in the context of \csa s.

The following terminology will be useful. Let
$\ \displaystyle{\ca = \lim_\to (A(G_k),\phi_k)}$\
be a limit algebra
as in Theorem 2.4 with diagonal subalgebra $\ \cd\ $.\  Refer to such an
algebra
as a {\em regular digraph limit algebra} and say that
$\ \cd$\  is a {\em regular canonical masa}, both of $\ \ca$\  and the
superalgebra
$\ \cb = C^*(\ca)$. In the self-adjoint context, $\ \cb = \ca\ $,\
for which we may asssume that each $\ G_k$\  is a union of
complete digraphs, it is known that a regular canonical masa
is independent of the presentation of $\ \ca\ $,\  in the following sense: if
$\ \cd$\  and $\ \cd^\prime$\  are two such masas in $\ \cb\ $,\
arising from different presentations of $\ \cb $,\ then
there is an approximately inner
automorphism $\ \alpha : \cb \to \cb$\  such that $\ \alpha(\cd) =
\cd^\prime$.\
This uniqueness theorem is due to Kreiger (see Renault \cite{ren})
and a direct proof
is given in \cite{scp-book}.
It would be very interesting to know if regular canonical  masas were
unique in this way in general (cf. Remark 1.8).
The following non-self-adjoint generalisation is straightforward.
\vspace{.3in}

\noindent  {\bf Theorem 4.1}\ \ \it Let $\ \ca =  A(G)
\otimes B$\  where $\ B$\
is an AF \csa \ and $\ A(G)$\  is a digraph algebra.
If $\ \cc$\  and $\ \cc^\prime$\  are regular canonical masas of $\ \ca$\
then there exists an approximately inner
automorphism $\ \alpha : \ca \to \ca$\  with $\ \alpha(\cc) = \cc^\prime\ $.\
\rm

\begin{mproof} We give a proof for the case when
$\ B$\ is an UHF \csa \ - the
setting for Theorem 4.5 -
and leave the reader to make the minor
changes necessary for the general case.

Assume that $\ G$\  is connected.
Let $\ \{h_{i,j}\}$\ be a partial matrix unit system for
$\ A(G)$.\
Note first that a masa
in $\ A(G) \ot B$\ is inner unitarily equivalent to one of the form
$\ h_{1,1} \ot C^{(1)} + \dots + h_{r,r} \ot C^{(r)}$\
where $\ r = |G|$\  and where each $\ C^{(k)}\  $ is a regular
canonical masa in $\ B$.\
We show now that we can further arrange that the masas $\ C^{(k)}\ $
coincide and are equal to a regular canonical masa, $\ C$\  say, in the
\csa \ $B$. Since $\ \cc^\prime $\ is similarly conjugate to a
masa of the form
$\ h_{1,1} \ot C^\prime + \dots + h_{r,r} \ot C^\prime $\
for some regular canonical masa $\ C^\prime$\ in $\ B$,\
the theorem follows readily from the self-adjoint case.

The masa  $\ \cc$\  can be described in the following way. There
is a matrix unit system
$\ \{e_{i,j}^{(k)}\}$\ for $\ \cb = C^*(\ca)$\ such that for each $\ k\ $
the finite system  $\ \{e_{i,j}^{(k)}\}$\  is a full matrix unit system for
a unital matrix
subalgebra $\ \cb_k \ $ of $\ \cb, \ $ and the following
properties hold :

(i) for fixed $\ k\ $each partial isometry $\ e^{k}_{p,q} \ $
is a sum of some of the matrix units of $\ \{e_{i,j}^{(k+1)}\},$\

(ii) the matrix algebra inclusions $\ \cb_k \subseteq
\cb_{k+1}\ $ are unital,

(iii) $\cc$\ is the closed span of a chain of masas
$\ \cc_k \subseteq \cb_k\ $ where $\ \cc_k = $\ span$\{e_{i,i}^{(k)}\},$

(iv) $\ca \cap \cb_k \ $ is spanned by some of the matrix
units of $\ \{e_{i,j}^{(k)}\},$\  including all the diagonal matrix units
$\ \{e_{i,i}^{(k)}\}.$\

Without loss of generality assume that $\ h_{j,j} \ot 1\ $
lies in $\ \cc_1 \ $ for each $\ j.\ $
Then each $\ h_{j,j} \ot 1\ $ is the sum of the
same number of minimal diagonal
matrix units in the set $\ \{e_{i,i}^{(1)}\}.\ $
It follows that there is a partial isometry $\ v\ $ in $\cb_1$\ which
is a sum of matrix units in the set  $\ \{e_{i,j}^{(1)}\}$\ and has
initial projection  $\ h_{2,2} \ot 1 \ $ and final projection
$\ h_{1,1} \ot 1. \ $ Necessarily $\ v = h_{1,2} \ot w \ $for some
partial isometry $\ w \ $ in $\ B.\ $ Since it is a sum of matrix units
it must normalise the masa $\ \cc \ $ and so
$\ v(h_{2,2} \ot C^{(2)})v^* = h_{1,1} \ot C^{(1)} \ $ and hence
$\ h_{1,1} \ot wC^{(2)}w^* = h_{1,1} \ot C^{(1)}.\ $
Using such elements $\ w\ $  construct a unitary operator in
the diagonal algebra $\ \sum h_{i,i} \ot B$ which conjugates
$\ \cc \ $ to a masa of the desired form.
\end{mproof}

As in the finite-dimensional setting, the following definition
is now well-defined and natural.
\vspace{.3in}

\noindent  {\bf Definition 4.2}\ \ For $\ \ca =  A(G) \otimes B$\
as above define
the {\em stable regular (partial isometry) homology} of $\ \ca$\  to
be the groups $\ H_n(\ca) = H_n(\ca;\cc),$\
for $\ n = 0, 1, 2,...,$ where  $\ \cc$\  is a
regular canonical masa of $\ \ca\ $.\
\vspace{.3in}

\noindent  {\bf Definition 4.3}\ \ Let $\ \ca, \ \ca^\prime$\  be regular
digraph
limit algebras. Then

(i)\ \ an algebra homomorphism $\ \ca \to \ca^\prime $\
is said to be {\em regular}
if there exist regular canonical masas $\ \cc \subseteq \ca,\ \
\cc^\prime \subseteq \ca^\prime$\  such that $\ \beta(\cc) \subseteq
\cc^\prime$\
and $\ \beta(N_\cc(\ca)) \subseteq N_{\cc^\prime}(\ca^\prime)\ $where $\
N_\cc(\ca)$\ is the partial isometry normaliser of $\ \cc$\ in $\ \ca$.

(ii)\ \ If $\ \ca^{\prime\prime} \subseteq \ca$\  then $\ \ca^{\prime\prime}$\
is said to be
a {\em regular subalgebra} if it is star-extendibly isomorphic
to a regular digraph limit algebra
and the inclusion map is regular.
\vspace{.3in}

The simplest regular subalgebras are the closed subalgebras
$\ \ca^{\prime\prime}$\
such that $\ \cc \subseteq \ca^{\prime\prime} \subseteq \ca$\  for some regular
canonical masa $\ \cc$\  of $\ \ca\ $.\  These may be thought of
as  the multiplicity one
subalgebras. They are automatically
regular digraph limit algebras, and they are described
in terms of subrelations of the
approximately finite semigroupoid $\ R(\ca;\cc)$\  associated with $\ \cc\ $.\
For details see Chapter 7 of \cite{scp-book}. On the other hand
the unital inclusion
$\ A(G) \otimes B \to A(G) \otimes B \otimes M_n$\  given
by $\ a \to a \otimes 1_n$\  is
a regular inclusion of finite multiplicity $\ n$\  in the sense that
the commutant of the range is isomorphic
to $\ M_n\ $.\

In general, in addition to the index of the inclusion,
we need $\ K$-theoretic data, stable homology data, and perhaps
other invariants in order to determine the conjugacy class.

Extending the earlier usage, say that an embedding
$\ \alpha : A(G) \otimes B \to A(G) \otimes B^\prime$\
 is {\em rigid }
if there is an identification $\ B^\prime = M_n \otimes B$\  such
that $\ \alpha(a) =
\phi(a) \otimes id_B$\  where $\ \phi$\  is a rigid embedding. The
 multiplicity
of $\ \alpha$\  is defined to be the multiplicity
of $\ \phi\ $.\

In fact
such embeddings and their multiplicities
may be characterised intrinsically, without reference to a postulated tensor
decomposition,
in terms of the fundamental topological binary relation $\ R(\ca^\pr)$\
for the
pair $\ (\ca^\pr, \cc^\pr).$
This fact is not needed below but we nevertheless indicate
this characterisation in the case of the 4-cycle $\ G = D_4$.\

Let $\ v_1, v_2,v_3, v_4$\ be the images of
$\ e_{1,3} \ot 1, e_{1,4} \ot 1, e_{2,3} \ot 1, e_{2,4} \ot 1 $\ under
the rigid embedding $\ \ga$.\
For each point $\ x$\ in the Gelfand space $\ M(\cc^\pr)$,\ which is
dominated by the initial projection of one of these
images, the partial isometries $\ v_i$\  determine
a subgraph of $\ R(\ca^\pr)$.\
A simple compactness argument shows that the
embedding is rigid if and only if each such subgraph
is a nondegenerate copy of $\ G$\ in the sense of being equivalent
to the canonical copies of $\ G$.\ (Of course, while all these copies
of $\ G$\ are equivalent
this equivalence need not respect the labellings
inherited from the partial isometries $\ v_i$.)\

We now generalise Proposition 3.4 and classify the rigid embeddings
betweeen  cycle algebras of the form $\ \ca = A(D_{2m}) \otimes B\ $\
where $\ B$\  is a UHF \csa.
The
following extra homological invariant is needed.
\vspace{.3in}

\noindent  {\bf Definition 4.4}\ \ The {\em scale} of the stable homology group
$\ H_1(\ca;\cc)$\  is the subset $\ \gS_1(\ca;\cc)$\  of elements
arising from cycles
associated with partial matrix unit
systems $\ \{e_{i,j}\}$\  with $\ e_{i,i} \in \ca $\  for all $\ i\ $.\
\vspace{.3in}

In the case of the cycle algebras $\ \ca = A(D_{2n}) \ot B$\  we may write
$\ \gS_1(\ca)$\  for the scale
and there is a natural identification

\[
(H_1(\ca),\gS_1(\ca)) = (K_0(B), [-1_B,1_B])
\]

\noindent where
$\ (K_0(B), 1_B) = (\IQ(n),1)$\  and $\ \IQ(n)$\  is the subgroup of $\ \IQ$\
associated
with the generalised integer $\ n$\  for $\ B\ $.\  Define the
{\it scale} of
$\ K_0\ca \op H_1\ca$\  to be the subset
of $\ \gS(\ca) \times \gS_1(\ca)$\
consisting of the pairs $\ ([p],\gs)$\ where $\ \gs \in \gS_1$\ arises from a
cycle associated with a partial matrix unit
system $\ \{e_{i,j}\}$\ with $\ [p] = [e_{i,i}]$.
\vspace{.3in}

\noindent  {\bf Theorem 4.5}\ \ \it Let
$\ \ca_1 = A(D_{2n}) \ot B$\  and
$\ \ca_2 = A(D_{2n}) \ot B^\prime$\
where $\ B$\  and $\ B^\pr$\  are UHF \csa s,
and let
$\ \ga_i : \ca_1 \to \ca_2,
i = 1,2,$\  be rigid embeddings.
Then $\ \ga_1$\  and $\ \ga_2$\  are inner unitarily
equivalent if and only if the following conditions hold.

(i) $\ \ga_1$\  and $\ \ga_2$\  have the same multiplicity.

(ii) $\ \ga_1$\  and $\ \ga_2$\  induce the same scaled group homomorphisms
from $\ K_0\ca_1 \op H_1\ca_1$\ to $\ K_0\ca_2 \op H_1\ca_2.$
\rm
\vspace{.3in}

\begin{mproof}
The necessity of the conditions is
straightforward.

For the converse we may assume, by replacing $\ \ga_1$\  and $\ \ga_2$\
by conjugate maps, that $\ \ga_i = \phi_i \ot id_B$\
where $\ B^\prime = M_m \ot B $\
for some integer $\ m\ $,\  which is greater than the multiplicities
of $\ \ga_1$\  and $\ \ga_2\ $,\ and
where each map $\ \phi_i : A(D_4) \to A(D_4) \ot M_m$\  is a rigid
 embedding. Thus, in view of Proposition 3.4 it remains to
show that the information of (i) and (ii) is sufficient
to determine $\ K_0\phi_i$\  and $\ H_1\phi_i\ $.\

Let $\ s$\  be the generalised integer for $\ B\ $.\  Then
$\ (K_0\ca, \gS(\ca))$\  is identifiable with the $\ 2n$-fold product

\[
(\ \IQ(s) \op \dots \op \ \IQ(s), [0,1]^{2n})
\]

\noindent and

\[
(H_1\ca,\gS_1(\ca)) = (\ZZ \ot_{\ZZ} \IQ(s),[-1,1]) = (\IQ(s),[-1,1]).
\]

\noindent There are similar identifications for $\
(K_0\ca^\prime,\gS(\ca^\prime))$\
with $\ ms$\  in place of $\ s$\  and under these identifications
it follows that $\ K_0\ga_i\ $,\
as a $\ 2n$\  by $\ 2n$\  matrix, is equal to $\ q_iK_0\phi_i\ $,\
where $\ q_i $\ is
equal to the inverse of the multiplicity of $\ \phi_i\ $.\
Furthermore, as a $\ 1$\  by $\ 1$\  matrix, $\ H_1\ga_i$\  is equal to
$\ q_iH_1\phi_i\ $.\
By the hypotheses it follows that $\ \phi_1$\  and $\ \phi_2$\  have coincident
$\ K_0 \op H_1$\  data, as desired.
\end{mproof}
\vspace{.3in}

\noindent  {\bf Remark 4.6}
\ \ As we have already mentioned
it would be desirable to generalise
Theorem 4.1 to general regular limits of digraph algebras.
The essential obstacle for this is already present in the case of algebraic
direct limits.
Suppose that $\ \ca$\ is such  a limit algebra with two
digraph subalgebra chains

\[
A_k \subseteq A_{k+1} \ \ \mbox{and}\ \
A_k^\pr \subseteq A_{k+1}^\pr  \ \ \ \
\]

\noindent for all $\ k = 1,2,\dots,$\ where $\ A_1, A_2,\dots$\ and
$\ A_1^\pr, A_2^\pr, \dots$\ are digraph algebras, with dense union, for
which the given inclusions are regular.
In particular it is possible to choose partial matrix
unit systems, in the usual sense,
for the chains $\ \{A_k\}$\ and $\ \{A_k^\pr\},$\ which in turn determine
regular canonical masas, $\ \cc $\ and $\ \cc^\pr$\ say,
spanned by the diagonal matrix units.
Choosing subsystems and relabelling we may assume furthermore
that $\ A_k \subseteq A_k^\pr \subseteq A_{k+1}$\ for all $\ k$.\
If these inclusions are regular then it can be shown that
$\ \cc$\ and $\ \cc^\pr$\ are conjugate by an approximately inner automorphism
of $\ \ca$.\ (In particular it follows that the
conjugacy class of $\ \cc$\ is determined
by the chain $\ \{A_k\}$ and is independent of the choice
of matrix unit system.)\ However examples can be constructed wherein these
inclusions are not regular.

\newpage

\section  {\bf Limit algebras}

The following discussion illustrates the use of
$\ K_0 \oplus H_1$-uniqueness in the identification of limit algebras.

Consider the system $\ \ca_1 \to \ca_2 \to \dots \ $\
consisting of $4$-cycle digraph algebras

\[
A(D_4) \op A(D_4) \to (A(D_4) \op A(D_4)) \ot M_{20} \to
(A(D_4) \op A(D_4)) \ot M_{20^2} \to \dots \ca.
\]

\noindent Assume furthermore that this is  a stationary
direct system in which each embedding is a fixed rigid embedding
similar to the
type mentioned before Definition 3.3. That is, the $k^{th}$
embedding of the system has the form
$\ \phi_k = \phi \ \ot \ $id$_{k-1} : \ca_1 \ot M_{10^{k-1}} \to (\ca_1 \ot
M_{10}) \ot M_{10^{k-1}} $\ where

\[
\phi = \left[ \begin{array}{cc}
\psi_1 & \psi_2 \\
\psi_3 & \psi_4
\end{array} \right]
\]

\noindent and where each partial embedding $\ \psi_i$\ is
a rigid  embedding of the form
$r_1\theta_1 + \dots + r_4\theta_4$.\ (The coefficients $\ r_k$\ depend on $\
i$.)\
Make the additional restriction that

\[
K_0\phi =
\left[
\begin{array}{cccccccc}
5&5&0&0&5&5&0&0\\
5&5&0&0&5&5&0&0\\
0&0&5&5&0&0&5&5\\
0&0&5&5&0&0&5&5\\
5&5&0&0&5&5&0&0\\
5&5&0&0&5&5&0&0\\
0&0&5&5&0&0&5&5\\
0&0&5&5&0&0&5&5
\end{array}
\right],
\]

\noindent so that

\[
H_0\phi =
\left[ \begin{array}{cc}
10 & 10\\10 & 10
\end{array} \right],
\]

\noindent and for convenience denote these matrices
by $\ T\ $ and $\ S$\ respectively.
With these assumptions the stationary limit algebra
$\ \ca $\ is determined by the $2 \times 2$ integral matrix
$\ X = H_1\phi.$\ Write $\ \ca_X$\ for the algebra.
For each of the partial embeddings $\ \psi$\ of $\ \phi$\
there are six possibilities. In the notation of Proposition 3.4 these are

\[
5\theta_1 + 5\theta_3, \ 4\theta_1 + 4\theta_3 + \theta_2 + \theta_4,
\ 3\theta_1 + 3\theta_3 + 2\theta_2 + 2\theta_4,
\]
\[
2\theta_1 + 2\theta_3 + 3\theta_2 + 3\theta_4,
\ \theta_1 + \theta_3 + 4\theta_2 + 4\theta_4,
\ 5\theta_2 + 5\theta_4.
\]

\noindent The induced homomorphisms on $\ H_1$\ are the
maps $\ \ZZ \to \ZZ$\ with entries

\[
10, 6, 2, -2, -6, -10,
\]

\noindent respectively. These numbers form the so called homology range
(in the terminology of \cite{scp-book} of a rigid embedding
for $\ K_0\psi$\ (and, by terminological extension, for $\ \psi$\ itself).
There are thus $\ 6^4$\ possibilities for the matrix $\ X$,\ and, a priori,
a great many possibilities for the limit algebras
$\ \ca_X$.\ Note that all of these algebras induces the same inclusion

\[
\ca_X  \cap \ca_X^* \to C^*(\ca_X).
\]

Let us focus on two of these algebras, namely

\[
\ca_{\tiny \left[ \begin{array}{cc}
10&6\\6&10\end{array}\right]  } \ \ \mbox{and}\ \
\ca_{\tiny \left[ \begin{array}{cc}
6&2\\2&6\end{array}\right]  }.
\]

\noindent This pair is of interest because, with respect to the natural masas,

\[
H_1(\ca_{\tiny \left[ \begin{array}{cc}
10&6\\6&10\end{array}\right]}) = H_1(\ca_{\tiny \left[ \begin{array}{cc}
6&2\\2&6\end{array}\right]}) = \ \ \IQ(2^{\infty}) \op \IQ(2^{\infty}).
\]

\noindent
Coincidence of this homology suggests that the two limit algebras
may be isomorphic, and indeed they are.

The method of proof in this rather typical stationary
example is to make use of the $\ K_0 \op H_1-$uniqueness property
to construct a commuting diagram
linking the two systems for the algebras.
\vspace{.3in}

\noindent {\bf Proposition 5.1}\ \ \it The 4-cycle limit algebras
\ $\ca_{\tiny \left[ \begin{array}{cc}
10&6\\6&10\end{array}\right]  } \ \ {and}\ \
\ca_{\tiny \left[ \begin{array}{cc}
6&2\\2&6\end{array}\right]  }$\
are star-extendibly isomorphic.
\rm
\vspace{.3in}

\begin{mproof}
Let $\ \displaystyle{\ca = \lim_\to (\ca_k, \phi_k)}, \
 \displaystyle{\ca^\prime = \lim_\to (\ca_k^\prime, \phi_k^\prime)}$\
be the respective systems for the algebras, as above,
and let $\ X$\ and $\ Y$\  be their respective 2 by 2 integral matrices.
Consider the commuting
diagram
\vspace{.3in}

\begin{center}
\setlength{\unitlength}{0.0125in}%
\begin{picture}(335,119)(30,660)
\thicklines
\put( 85,680){\vector( 1, 0){ 40}}
\put(180,680){\vector( 1, 0){ 40}}
\put(305,680){\vector( 1, 0){ 40}}
\put( 85,710){\vector( 4, 3){ 40}}
\put(180,760){\vector( 3,-1){165}}
\put( 50,745){\vector( 0,-1){ 40}}
\put( 40,760){\makebox(0,0)[lb]{\raisebox{0pt}[0pt][0pt]{$H_1(\ca_1)$}}}
\put( 40,680){\makebox(0,0)[lb]{\raisebox{0pt}[0pt][0pt]{$H_1(\ca_1^\pr)$}}}
\put(135,680){\makebox(0,0)[lb]{\raisebox{0pt}[0pt][0pt]{$H_1(\ca_2^\pr)$}}}
\put( 85,760){\vector( 1, 0){ 40}}
\put(250,680){\makebox(0,0)[lb]{\raisebox{0pt}[0pt][0pt]{\twlrm . . . . . . }}}
\put( 30,725){\makebox(0,0)[lb]{\raisebox{0pt}[0pt][0pt]{\twlrm id}}}
\put(365,680){\makebox(0,0)[lb]{\raisebox{0pt}[0pt][0pt]{$H_1(\ca_{2+j}^\pr)$}}}
\put(135,760){\makebox(0,0)[lb]{\raisebox{0pt}[0pt][0pt]{$H_1(\ca_2)$}}}
\put( 95,770){\makebox(0,0)[lb]{\raisebox{0pt}[0pt][0pt]{$X$}}}
\put( 90,660){\makebox(0,0)[lb]{\raisebox{0pt}[0pt][0pt]{$Y$}}}
\put(185,660){\makebox(0,0)[lb]{\raisebox{0pt}[0pt][0pt]{$Y$}}}
\put(320,660){\makebox(0,0)[lb]{\raisebox{0pt}[0pt][0pt]{$Y$}}}
\put( 85,725){\makebox(0,0)[lb]{\raisebox{0pt}[0pt][0pt]{$U_1$}}}
\put(280,730){\makebox(0,0)[lb]{\raisebox{0pt}[0pt][0pt]{$V_1$}}}
\end{picture}
\end{center}

\vspace{.3in}

\noindent where $\ U_1 = X.$\  We wish to choose $\ j$\ large
enough so that the matrix $\ V_1 = Y^{1+j}U_1^{-1}$\ is an integral matrix
belonging  to the homology
range of the map

\begin{center}
\setlength{\unitlength}{0.0125in}%
\begin{picture}(215,24)(120,755)
\put(120,755){\makebox(0,0)[lb]{\raisebox{0pt}[0pt][0pt]{$K_0\ca_2$}}}
\thicklines
\put(160,760){\vector( 1, 0){160}}
\put(335,755){\makebox(0,0)[lb]{\raisebox{0pt}[0pt][0pt]
{$K_0\ca_{2+j}^\prime$}}}
\put(220,770){\makebox(0,0)[lb]{\raisebox{0pt}[0pt][0pt]{$T^j$}}}
\end{picture}
\end{center}

\noindent Note that the homology range can be easily calculated from the matrix
$\ S^j$.\
In fact $\ j = 2 $\ is the first index for which this occurs, with

\[
V_1  = \left[\begin{array}{cc}
24&8\\8&24\end{array}\right]
\ \ \ \ \
S^2 = \left[\begin{array}{cc}
200&200\\200&200\end{array}\right]
\]

\noindent We can now simultaneously lift $\ U_1$\
and $\ T$\ to a rigid embedding
$\ \gb : \ca_1^\prime \to \ca_2 $\ and we can lift $\ V_1$\ and $\ T^2$\ to
a rigid embedding
$\ga_1 : \ca_2 \to \ca_4^\prime$.\ Furthermore since

\[
K_0 \op H_1(\ga_1 \circ \gb_1) =
K_0 \op H_1(\phi_3^\pr \circ \phi_2^\pr \circ \phi_1^\pr)
\]

\noindent we may apply Proposition 3.4 and replace $\ \ga_1$\ by an
inner conjugate map so that

\[
\ga_1 \circ \gb_1 = \phi_3^\pr \circ \phi_2^\pr \circ \phi_1^\pr
\]

Consider next the diagram
\vspace{.3in}

\begin{center}
\setlength{\unitlength}{0.0125in}%
\begin{picture}(240,104)(125,600)
\thicklines
\put(180,680){\vector( 1, 0){ 40}}
\put(180,680){\vector( 1, 0){ 40}}
\put(180,680){\vector( 1, 0){ 40}}
\put(150,665){\vector( 0,-1){ 40}}
\put(180,605){\vector( 3, 1){165}}
\put(135,680){\makebox(0,0)[lb]{\raisebox{0pt}[0pt][0pt]{$H_1\ca_2$}}}
\put(305,680){\vector( 1, 0){ 40}}
\put(250,680){\makebox(0,0)[lb]{\raisebox{0pt}[0pt][0pt]{\twlrm . . . . . . }}}
\put(125,640){\makebox(0,0)[lb]{\raisebox{0pt}[0pt][0pt]{$ V_1$}}}
\put(365,680){\makebox(0,0)[lb]{\raisebox{0pt}[0pt][0pt]{$H_1\ca_{2+k}$}}}
\put(135,600){\makebox(0,0)[lb]{\raisebox{0pt}[0pt][0pt]{$ H_1\ca_4^\pr$}}}
\put(280,615){\makebox(0,0)[lb]{\raisebox{0pt}[0pt][0pt]{$U_2$}}}
\put(195,695){\makebox(0,0)[lb]{\raisebox{0pt}[0pt][0pt]{$X$}}}
\put(320,695){\makebox(0,0)[lb]{\raisebox{0pt}[0pt][0pt]{$X$}}}
\end{picture}
\end{center}
\vspace{.3in}

We wish to choose $\ k$\ large enough so that the
matrix $\ U_2 = X^kV_1^{-1}$\ is an
integral matrix lying in the homology range of
$\ T^k$.\  It is clear that such a $\ k$\ exists for the following
two reasons.

(i) the entries of $\ T^k$\ will eventually exceed in modulus
the corresponding entries of $X^kV_1^{-1}.$

(ii) all entries of $\ T^k$\ and $X^kV_1^{-1}$\ are congruent to zero mod 4
for sufficiently large $\ k$.

\noindent In fact the first value for which (i) and (ii) hold
is $\ k = 4$ giving

\[
X^4V_1^{-1} = \left[\begin{array}{cc}
1032&1016\\1016&1032\end{array}\right],
\ \ \ \ \ \ \
S^4 = \left[\begin{array}{cc}
80000&80000\\80000&80000\end{array}\right].
\]

\noindent As before we can lift $\ U_2$\ to a rigid homomorphism
$\ \gb_2$\ in such a way that we obtain a commuting triangle  so that
$\ \gb_2 \circ \ga_1 = \phi_5 \circ \phi_4 \circ \phi_3 \circ \phi_2 $.\
It is clear that the requirements of (i) and (ii) can always be met
at further stages
in the construction of the commuting diagram. In this way we
obtain the desired  commuting diagram

\vspace{.3in}

\begin{center}
\setlength{\unitlength}{0.0125in}%
\begin{picture}(449,104)(80,680)
\thicklines
\put(100,740){\vector( 0,-1){ 40}}
\put(240,740){\vector( 0,-1){ 40}}
\put(140,705){\vector( 3, 2){ 57.692}}
\put(140,765){\vector( 1, 0){ 60}}
\put(140,685){\vector( 1, 0){ 60}}
\put(280,765){\vector( 1, 0){ 60}}
\put(280,685){\vector( 1, 0){ 60}}
\put(280,705){\vector( 3, 2){ 57.692}}
\put(425,740){\vector( 0,-1){ 40}}
\put( 90,760){\makebox(0,0)[lb]{\raisebox{0pt}[0pt][0pt]{$\ca_1$}}}
\put( 90,680){\makebox(0,0)[lb]{\raisebox{0pt}[0pt][0pt]{$\ca_{n_1}^\prime$}}}
\put(235,760){\makebox(0,0)[lb]{\raisebox{0pt}[0pt][0pt]{$\ca_{m_1}$}}}
\put(235,680){\makebox(0,0)[lb]{\raisebox{0pt}[0pt][0pt]{$\ca_{n_2}^\prime$}}}
\put( 80,715){\makebox(0,0)[lb]{\raisebox{0pt}[0pt][0pt]{$\ga_1$}}}
\put(220,715){\makebox(0,0)[lb]{\raisebox{0pt}[0pt][0pt]{$\ga_2$}}}
\put(375,725){\makebox(0,0)[lb]{\raisebox{0pt}[0pt][0pt]{$\dots$}}}
\put(415,760){\makebox(0,0)[lb]{\raisebox{0pt}[0pt][0pt]{$ \ca$}}}
\put(415,680){\makebox(0,0)[lb]{\raisebox{0pt}[0pt][0pt]{$ \ca^\prime$}}}
\put(405,720){\makebox(0,0)[lb]{\raisebox{0pt}[0pt][0pt]{$ \ga$}}}
\end{picture}
\end{center}
\vspace{.3in}

\end{mproof}

The reader may notice that the stationary case above presents no
difficulties with regard to the harmonisation of the homology coupling
invariants given in Chapter 11 of \cite{scp-book}. Addressing this issue
is just one of the tasks necessary for a complete classification of
rigid embedding limits of digraph algebras.

Using the method of the last proof one can obtain the following more
general theorem.
\vspace{.3in}

{\bf Theorem 5.2} \ \ \it Let $\ca_X$\ and $\ca_Y$\ be limit algebras,
as above, associated with a pair of 2 by 2 matrices whose entries
lie in the set \ $\{10, 6, 2, -2, -6, -10\}$.\ If the (diagonal masa)
homology groups $\ H_1(\ca_X)$\ and $\ H_1(\ca_Y)$\ are isomorphic
then $\ \ca_X$\ and $\ca_Y$\ are star-extendibly isomorphic
operator algebras.
Furthermore, for the algebras $\ca_X$ \ with
$X = {\tiny \left[ \begin{array}{cc}
a&b\\b&a\end{array}\right]  }$
there are at most five
isomorphism classes corresponding to the groups
$\ \IQ(2^\infty),
\ \IQ(6^\infty),
\ \IQ(10^\infty),
\ \IQ(2^\infty) \op  \IQ(2^\infty),
\ \IQ(2^\infty) \op  \IQ(6^\infty).
$\
\rm

\newpage

\begin{center}
\large
\bf Appendix 1 \rm

The Coefficient Matrix for the rotation embeddings of $Cu$
\end{center}
\vspace{.3in}

\noindent Label $\ Cu$ \  in the following manner, where
the receiving vertices are labelled \ 1,2,3,4.

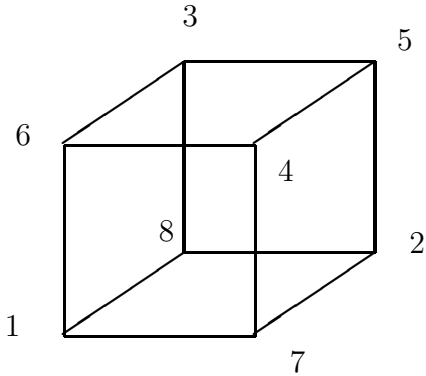
\begin{center}
\setlength{\unitlength}{0.0125in}%
\begin{picture}(165,159)(135,580)
\thicklines
\put(210,715){\line( 1, 0){ 80}}
\put(290,715){\line( 0, 1){  0}}
\put(290,715){\line( 0,-1){ 80}}
\put(290,635){\line( 0, 1){  0}}
\put(290,635){\line(-1, 0){ 80}}
\put(210,635){\line( 0, 1){  0}}
\put(210,635){\line( 0, 1){ 80}}
\put(210,715){\line(-3,-2){ 50.769}}
\put(290,715){\line(-3,-2){ 50.769}}
\put(290,635){\line(-3,-2){ 50.769}}
\put(210,635){\line(-3,-2){ 50.769}}
\put(135,600){\makebox(0,0)[lb]{\raisebox{0pt}[0pt][0pt]{\twlrm 1}}}
\put(160,680){\line( 1, 0){ 80}}
\put(240,680){\line( 0, 1){  0}}
\put(240,680){\line( 0,-1){ 80}}
\put(240,600){\line( 0, 1){  0}}
\put(240,600){\line(-1, 0){ 80}}
\put(160,600){\line( 0, 1){  0}}
\put(160,600){\line( 0, 1){ 80}}
\put(300,635){\makebox(0,0)[lb]{\raisebox{0pt}[0pt][0pt]{ 2}}}
\put(195,640){\makebox(0,0)[lb]{\raisebox{0pt}[0pt][0pt]{ 8}}}
\put(205,730){\makebox(0,0)[lb]{\raisebox{0pt}[0pt][0pt]{ 3}}}
\put(245,665){\makebox(0,0)[lb]{\raisebox{0pt}[0pt][0pt]{ 4}}}
\put(295,720){\makebox(0,0)[lb]{\raisebox{0pt}[0pt][0pt]{ 5}}}
\put(135,680){\makebox(0,0)[lb]{\raisebox{0pt}[0pt][0pt]{ 6}}}
\put(250,585){\makebox(0,0)[lb]{\raisebox{0pt}[0pt][0pt]{ 7}}}
\end{picture}
\end{center}

\noindent The $\ K_0$ maps of the 12 multiplicity one embeddings
associated with the \ 12\  rotations of $\ Cu$\ are given by

\[
T_1 = \left [\begin {array}{cccc} 1&0&0&0\\0&0&1&0\\0&0&0&1\\0&1&0&0
\end {array}\right ],\ \
T_2 = \left [\begin {array}{cccc} 0&0&1&0\\0&1&0&0\\0&0&0&1\\1&0&0&0
\end {array}\right ],\ \
T_3 = \left [\begin {array}{cccc} 0&0&0&1\\0&0&1&0\\0&1&0&0\\1&0&0&0
\end {array}\right ],\ \\
\]
\[T_4 = \left [\begin {array}{cccc} 0&1&0&0\\1&0&0&0\\0&0&0&1\\0&0&1&0
\end {array}\right ],\ \
T_5 = \left [\begin {array}{cccc} 0&0&1&0\\0&0&0&1\\1&0&0&0\\0&1&0&0
\end {array}\right ],\ \
T_6 = \left [\begin {array}{cccc} 1&0&0&0\\0&0&0&1\\0&1&0&0\\0&0&1&0
\end {array}\right ],
\]\[
T_7 = \left [\begin {array}{cccc} 0&0&0&1\\1&0&0&0\\0&0&1&0\\0&1&0&0
\end {array}\right ],\ \
T_8 = \left [\begin {array}{cccc} 0&0&0&1\\0&1&0&0\\1&0&0&0\\0&0&1&0
\end {array}\right ],\ \
T_9 = \left [\begin {array}{cccc} 0&1&0&0\\0&0&0&1\\0&0&1&0\\1&0&0&0
\end {array}\right ],
\]\[
T_{10} = \left [\begin {array}{cccc} 0&0&1&0\\1&0&0&0\\0&1&0&0\\0&0&0&1
\end {array}\right ],\ \
T_{11} = \left [\begin {array}{cccc} 1&0&0&0\\0&1&0&0\\0&0&1&0\\0&0&0&1
\end {array}\right ],\ \
T_{12} = \left [\begin {array}{cccc} 0&1&0&0\\0&0&1&0\\1&0&0&0\\0&0&0&1
\end {array}\right ].
\]
\rm
\vspace{.3in}
\rm

\noindent Consider the basis of $\ H_1(A(Cu)) = \ZZ^5$\ given by the cycles

\[
<3,5> + <5,4> + <4,6> + <6,3>,
\] \[
<6,2> + <2,8> + <8,1> + <1,6>,
\] \[
<3,5> + <5,2> + <2,8> + <8,3>,
\] \[
<5,4> + <4,7> + <7,2> + <2,5>,
\] \[
<4,6> + <6,1> + <1,7> + <7,4>.
\]

\noindent Then the following matrices represent the corresponding $\ H_1$\
maps
of the 12 multiplicity one rotation embeddings.

\[
S_1 = \left [\begin {array}{ccccc} 0&-1&0&1&0\\0&-1&0&0&1\\1&-1&0&0&0\\0&-1&
1&0&0\\0&-1&0&0&0\end {array}\right ],\ \
S_2 = \left [\begin {array}{ccccc} 0&0&-1&0&1\\1&0&-1&0&0\\0&0&-1&1&0\\0&0&-
1&0&0\\0&1&-1&0&0\end {array}\right ],\ \
S_3 = \left [\begin {array}{ccccc} -1&0&0&0&0\\-1&0&0&1&0\\-1&0&1&0&0\\-1&1&0
&0&0\\-1&0&0&0&1\end {array}\right ],
\]\[
S_4 = \left [\begin {array}{ccccc} 1&0&0&0&0\\0&0&0&1&0\\0&0&0&0&1\\0&1&0&0&0
\\0&0&1&0&0\end {array}\right ],\ \
S_5 = \left [\begin {array}{ccccc} -1&0&0&0&0\\-1&1&0&0&0\\-1&0&0&0&1\\-1&0&0
&1&0\\-1&0&1&0&0\end {array}\right ],\ \
S_6 = \left [\begin {array}{ccccc} 0&0&1&0&-1\\0&0&0&0&-1\\0&0&0&1&-1\\1&0&0
&0&-1\\0&1&0&0&-1\end {array}\right ],
\]\[
S_7 = \left [\begin {array}{ccccc} 0&0&1&0&-1\\1&0&0&0&-1\\0&1&0&0&-1\\0&0&0
&0&-1\\0&0&0&1&-1\end {array}\right ],\ \
S_8 = \left [\begin {array}{ccccc} 0&0&1&0&-1\\1&0&0&0&-1\\0&1&0&0&-1\\0&0&0
&0&-1\\0&0&0&1&-1\end {array}\right ],\ \
S_9 = \left [\begin {array}{ccccc} 0&1&0&-1&0\\0&0&1&-1&0\\1&0&0&-1&0\\0&0&0
&-1&1\\0&0&0&-1&0\end {array}\right ]
\]\[
S_{10} = \left [\begin {array}{ccccc}
0&-1&0&1&0\\0&-1&1&0&0\\0&-1&0&0&0\\0&-1&0
&0&1\\1&-1&0&0&0\end {array}\right ],\ \
S_{11} = \left [\begin {array}{ccccc}
1&0&0&0&0\\0&1&0&0&0\\0&0&1&0&0\\0&0&0&1&0
\\0&0&0&0&1\end {array}\right ],\ \
S_{12} = \left [\begin {array}{ccccc} 0&0&-1&0&1\\0&0&-1&0&0\\0&1&-1&0&0\\1&0&-
1&0&0\\0&0&-1&1&0\end {array}\right ].
\]

\vspace{.3in}
\rm

\noindent The coefficient matrix in section 3 arises from the 41 equations
in the multiplicities $\ r_1,\dots,r_{12}$\  coming from the matrix equations

\[
r_1T_1 + \dots + r_{12}T_{12} = T
\]\[
r_1S_1 + \dots + r_{12}S_{12} = S.
\]

\newpage

\begin{center}
\large
\bf Appendix 2 \rm
\vspace{.3in}

Coefficient Matrix for the Rigid Embeddings of $Cu$

\tiny
\[
\left [\begin {array}{cccccccccccccccccccccccc} 1&0&0&0&0&1&0&0&0&0&1&0
&1&0&0&0&0&1&0&0&0&0&1&0\\0&0&0&1&0&0&0&0&1&0&0&1&0&0&0&1&0&0&0&0&1&0&0
&1\\0&1&0&0&1&0&0&0&0&1&0&0&0&1&0&0&1&0&0&0&0&1&0&0\\0&0&1&0&0&0&1&1&0
&0&0&0&0&0&1&0&0&0&1&1&0&0&0&0\\0&0&0&1&0&0&1&0&0&1&0&0&0&0&0&1&0&0&1&0
&0&1&0&0\\0&1&0&0&0&0&0&1&0&0&1&0&0&1&0&0&0&0&0&1&0&0&1&0\\1&0&1&0&0&0
&0&0&0&0&0&1&1&0&1&0&0&0&0&0&0&0&0&1\\0&0&0&0&1&1&0&0&1&0&0&0&0&0&0&0&
1&1&0&0&1&0&0&0\\0&0&0&0&1&0&0&1&0&0&0&1&0&1&1&0&0&0&0&0&1&0&0&0\\0&0&
1&0&0&1&0&0&0&1&0&0&1&0&0&0&1&0&1&0&0&0&0&0\\0&0&0&0&0&0&1&0&1&0&1&0&0
&0&0&1&0&1&0&1&0&0&0&0\\1&1&0&1&0&0&0&0&0&0&0&0&0&0&0&0&0&0&0&0&0&1&1&
1\\0&1&1&0&0&0&0&0&1&0&0&0&0&0&0&0&1&0&0&1&0&0&0&1\\1&0&0&0&1&0&1&0&0&0
&0&0&0&0&1&0&0&1&0&0&0&1&0&0\\0&0&0&1&0&1&0&1&0&0&0&0&0&0&0&0&0&0&1&0&
1&0&1&0\\0&0&0&0&0&0&0&0&0&1&1&1&1&1&0&1&0&0&0&0&0&0&0&0\\0&0&-1&1&-1&0
&0&0&0&0&1&0&0&0&1&-1&1&0&0&0&0&0&-1&0\\-1&0&0&0&0&0&0&1&1&-1&0&0&1&0&0
&0&0&0&0&-1&-1&1&0&0\\0&-1&0&0&0&1&1&0&0&0&0&-1&0&1&0&0&0&-1&-1&0&0&0&0
&1\\1&0&0&0&0&0&0&-1&-1&1&0&0&-1&0&0&0&0&0&0&1&1&-1&0&0\\0&1&0&0&0&-1&
-1&0&0&0&0&1&0&-1&0&0&0&1&1&0&0&0&0&-1\\0&1&-1&0&-1&0&1&0&0&0&0&0&0&0&
1&0&1&0&0&-1&0&-1&0&0\\-1&0&0&0&1&0&0&0&0&-1&1&0&1&-1&0&0&0&-1&0&0&0&1
&0&0\\0&-1&0&0&0&0&0&0&1&1&0&-1&0&1&0&-1&-1&0&0&0&0&0&0&1\\0&0&1&1&0&0
&0&-1&-1&0&0&0&0&0&0&0&0&0&-1&1&1&0&0&-1\\1&0&0&0&0&-1&-1&1&0&0&0&0&0&0
&-1&0&0&1&1&0&0&0&-1&0\\1&0&-1&0&-1&0&0&0&1&0&0&0&0&0&1&0&1&-1&0&0&0&0
&0&-1\\-1&0&0&0&0&0&1&0&0&-1&0&1&1&0&-1&-1&0&0&0&0&0&1&0&0\\0&-1&1&0&0
&0&0&0&0&0&1&-1&-1&1&0&0&0&0&0&-1&0&0&0&1\\0&1&0&0&0&1&0&-1&-1&0&0&0&0
&0&0&0&-1&0&0&1&1&0&-1&0\\0&0&0&1&1&-1&-1&0&0&0&0&0&0&0&0&0&0&1&1&0&-1
&-1&0&0\\0&0&-1&0&-1&1&0&0&0&0&0&1&-1&0&1&0&1&0&0&0&-1&0&0&0\\-1&0&1&1
&0&0&0&0&0&-1&0&0&1&0&0&0&0&0&-1&0&0&1&0&-1\\1&-1&0&0&0&0&0&1&0&0&0&-1
&0&1&-1&0&0&0&0&0&0&0&-1&1\\0&0&0&0&1&0&0&-1&-1&0&1&0&0&-1&0&0&0&-1&0&
1&1&0&0&0\\0&0&0&0&0&-1&-1&0&1&1&0&0&0&0&0&-1&-1&1&1&0&0&0&0&0\\0&0&-1
&0&-1&0&0&1&0&1&0&0&0&-1&1&0&1&0&-1&0&0&0&0&0\\-1&1&0&0&0&1&0&0&0&-1&0
&0&1&0&0&0&-1&0&0&0&0&1&-1&0\\0&-1&0&1&1&0&0&0&0&0&0&-1&0&1&0&0&0&0&0&0
&-1&-1&0&1\\0&0&0&0&0&0&1&-1&-1&0&0&1&0&0&-1&-1&0&0&0&1&1&0&0&0\\0&0&1
&0&0&-1&-1&0&0&0&1&0&-1&0&0&0&0&1&1&-1&0&0&0&0\end {array}\right ]
\]
\vspace{.3in}
\rm

\large The rank of this matrix is $\ 23$.
\end{center}

\end{document}